\documentclass[a4paper, 11pt]{article}

\usepackage{amsmath}
\usepackage{amsfonts}
\usepackage{amssymb}
\usepackage{amsthm}
\usepackage{bm}
\usepackage{indentfirst}
\usepackage{graphicx}
\usepackage{subfigure}
\usepackage{array}
\usepackage{color}
\usepackage{fullpage}
\usepackage{algorithm}
\usepackage{algorithmic}
\usepackage{url}
\usepackage{authblk}

\newcommand{\mo}{\mathcal{O}}
\DeclareMathAlphabet{\mathsfsl}{OT1}{cmss}{m}{sl}

\newcommand{\tensor}[1]{\mathbf{#1}}

\newcommand{\dif}{\mathrm{d}}

\newcommand{\me}{\mathrm{e}}
\newcommand{\ba}{\bm\alpha}
\newcommand{\bx}{\bm\xi}
\newcommand{\ve}{\varepsilon}
 
\newcommand{\mexp}[1]{\mathbb{E}\left\{{#1}\right\}}

\newcommand{\mr}{\mathbb{R}}
\newcommand{\mn}{\mathbb{N}}
\newcommand{\hc}{\hat{\bm c}} 
\newcommand{\trans}{\mathsf{T}}
\DeclareMathOperator*{\argmin}{arg\,min}
\newcommand{\PreserveBackslash}[1]{\let\temp=\\#1\let\\=\temp}
\newcolumntype{C}[1]{>{\PreserveBackslash\centering}p{#1}}
\newcolumntype{R}[1]{>{\PreserveBackslash\raggedleft}p{#1}}
\newcolumntype{L}[1]{>{\PreserveBackslash\raggedright}p{#1}}

\newcommand\Def{\stackrel{\textrm{def}}{=}}
\DeclareMathOperator{\sech}{sech}

\title{Sliced-Inverse-Regression-Aided Rotated Compressive Sensing Method for
Uncertainty Quantification}

\author{Xiu Yang\footnote{xiu.yang@pnnl.gov}}
\author{Weixuan Li\footnote{weixuan.li@pnnl.gov}}
\author{Alexandre Tartakovsky\footnote{alexandre.tartakovsky@pnnl.gov}}
\affil{Advanced Computing, Mathematics and Data Division, Pacific
Northwest National Laboratory, Richland, WA 99352}

\usepackage{amsopn}

\begin{document}
\maketitle

\begin{abstract}
Compressive-sensing-based uncertainty quantification methods have become a
powerful tool for problems with limited data. In this work, we use the sliced
inverse regression (SIR) method to provide an initial guess for the alternating
direction method, which is used to enhance sparsity of the Hermite polynomial 
expansion of stochastic quantity of interest. 
The sparsity improvement increases both the efficiency and accuracy of the 
compressive-sensing-based uncertainty quantification method. We demonstrate that
the initial guess from SIR is suitable for cases when the available data 
are limited (Algorithm \ref{algo:cs_sir1}). We also propose another algorithm
(Algorithm \ref{algo:cs_sir2}) that performs dimension reduction first with 
SIR. Then it constructs a Hermite polynomial expansion of the reduced model. 
This method affords the ability to approximate the statistics accurately 
with even less available data. Both methods are non-intrusive and require 
no \emph{a priori} information of the sparsity of the system. The
effectiveness of these two methods (Algorithms \ref{algo:cs_sir1} and
\ref{algo:cs_sir2}) are demonstrated using problems with up to $500$ random 
dimensions. 

\end{abstract}


\section{Introduction}
\label{sec:intro}

Surrogate model is a powerful tool in studying uncertainty quantification
(UQ). For example, spectral-method-based surrogate models, including the 
polynomial chaos expansion (PCE) \cite{GhanemS91} and generalized polynomial 
chaos (gPC) \cite{XiuK02} methods, are widely used for UQ in engineering and 
computational sciences. In the gPC and PCE methods, a quantity of interest (QoI) $u$
(e.g, velocity, temperature, etc.) depends on $d$-dimensional
($d<\infty$) random variables $\bm\xi=(\xi_1,\xi_2,\cdots,\xi_d)^\trans$, 
which are used to represent stochastic initial and boundary conditions 
or other unknown properties, can be approximated as 
\begin{equation}\label{eq:gpc}
u(\bx) = \sum_{n=1}^Nc_n\psi_n(\bx) + \ve(\bx),
\end{equation}
where $\ve$ is the truncation error; $N$ is a positive integer; $c_n$ are
coefficients; and $\psi_n$ are multivariate orthonormal polynomials satisfying
\begin{equation}
  \mexp{\psi_i(\bx)\psi_j(\bx)}=
  \int_{\mr^d} \psi_i(\bm x)\psi_j(\bm x)\rho_{_{\bx}}(\bm x)\dif\bm x = \delta_{i,j},
\end{equation}
where $\rho_{_{\bx}}(\bm x)$ is the probability density function (PDF) of $\bm\xi$
and $\delta_{ij}$ is the Kronecker delta function. 
Here $\bm\xi$ is defined on the probability space $(\Omega, \mathcal{F}, P)$,
where $\Omega$ is the abstract set of elementary events, $\mathcal{F}$ is a
$\sigma$-algebra of subsets of $\Omega$ and $P$ is the probability measure
on $\mathcal{F}$. The QoI $u$ is defined on the Hilbert space 
$\mathcal{H}=L^2(\Omega, \mathcal{F}, P)$ that consists of real-valued
random variables defined on $(\Omega, \mathcal{F}, P)$ with finite second moment
and is equipped with a inner product 
$(u,v)_{L^2} = \int_{\Omega}uv\dif P$, for $u, v\in \mathcal{H}$. 
For example, when $\{\xi_i\}_{i=1}^d$ are independent and identically 
distributed (i.i.d.) Gaussian random variables, i.e., 
$\bm\xi\sim \mathcal{N}(\bm 0, \tensor I)$, $P$ is the Gaussian
measure, and $\psi_n$ are normalized multi-variate Hermite polynomials. For this
case, systematic studies of convergence of PCE and 
gPC \cite{CameronM1947, ErnstMS12} indicate that $\sum_{n=1}^Nc_n\psi_n(\bx)$ 
converges to $u(\bm\xi)$ in $L^2$ as $N\rightarrow\infty$. For the convergence
of more general cases, we refer the interested readers to \cite{ErnstMS12}.

Both intrusive and non-intrusive methods 
\cite{GhanemS91,XiuK02,TatangPPM97, XiuH05,FooWK08,BabuskaNT10} are extensively 
used to compute the gPC coefficients $\bm c=(c_1,c_2,\cdots,c_N)^\trans$. 
Non-intrusive methods are more useful when the model used to obtain $u$ is
especially complex. These methods utilize training sets $\{(\bm\xi^q, u^q)\}_{q=1}^M$
to approximate coefficients $\bm c$. Here, $\bm\xi^q$ are samples of input based on
$\rho_{_{\bx}}$, and $u^q$ are corresponding samples of the output $u^q=u(\bx^q)$
obtained from the computational model. 
In many applications, it can be very costly to obtain $u^q$. Because of this, 
it often is $M<N$ or even $M\ll N$, making the following linear system 
underdetermined:
\begin{equation}\label{eq:cs_eq}
\tensor\Psi \bm c = \bm u - \bm\ve,
\end{equation}
where $\bm u=(u^1,u^2,\cdots,u^M)^\trans$ is the vector of output samples, 
$\tensor\Psi$ is an $M\times N$ matrix with $\Psi_{ij}=\psi_j(\bx^i)$
(where $j=1,\cdots, N$ and $i=1,\cdots, M$), and 
$\bm\ve=(\ve^1,\ve^2,\cdots,\ve^M)^\trans$ is a vector of error samples with 
$\ve^q=\ve(\bx^q)$ (where $q=1,\cdots, M$). 
The compressive sensing method has been shown to be effective at solving the 
underdetermined Eq.~\eqref{eq:cs_eq} when $\bm c$ is sparse 
\cite{CandesRT06,DonohoET06,Candes08,BrucksteinDE09}, i.e., solving
\begin{equation*}
(P_{1,\epsilon}):\qquad \argmin_{\hat{\bm c}}\Vert\hat{\bm c}\Vert_1, 
\text{~~subject to~~} \Vert\tensor\Psi\hat{\bm c}-\bm u\Vert_2\leq\epsilon,
\end{equation*}
to approximate $\bm c$ in Eq.~\eqref{eq:cs_eq} with $\hat{\bm c}$ (see 
Section 2.2). It has been used to solve UQ problems in various settings
\cite{DoostanO11,YanGX12,YangK13,LeiYZLB15,PengHD14,HamptonD15,JakemanES14,
SargsyanSNDRT14,PengHD15}. 

Several approaches have been developed to enhance the efficiency of 
solving Eq.~\eqref{eq:cs_eq} in UQ applications, including weighted/re-weighted $\ell_1$ 
minimization, which assigns a weight to each $c_n$ and solves a weighted $\ell_1$ 
minimization problem to enhance the sparsity \cite{CandesWB08,YangK13,PengHD14,RauhutW15};
smart sampling strategies to better the property of $\tensor\Psi$
\cite{RauhutW12,HamptonD15}; and adaptive basis selection 
to reduce the number of unknowns \cite{JakemanES14}. 

In \cite{LeiYZLB15,YangLBL16, YangWL17}, an approach to enhance the sparsity of
$\bm c$ through the rotation of the random vector $\bx$ has been proposed. 
This method aims to find a rotation $g:\mr^d\mapsto\mr^d$ that maps $\bx$ to a new 
set of random variables $\bm\eta=(\eta_1, \eta_2, \cdots, \eta_d)^\trans$ 
as $\bm\eta=g(\bm\xi)=\tensor A\bm\xi$ (where $\tensor A\tensor A^\trans=\tensor I$) such that
the gPC expansion of $u$ with respect to $\bm\eta$ is sparser. Specifically,
\begin{equation}\label{eq:gpc2}
u(\bx) \approx \sum_{n=1}^N c_n\psi_n(\bx)=\sum_{n=1}^N \tilde c_n \psi_n(\bm\eta(\bx))\approx u(\bm\eta(\bx)),
\end{equation}
and $\tilde{\bm c}=(\tilde c_1, \tilde c_2, \cdots, \tilde c_N)^\trans$ is
sparser than $\bm c$. Hence, $\tilde{\bm c}$ can be approximated more accurately
using the compressive sensing method. Subsequently, the enhancement of the sparsity
enables the compressive sensing algorithm to obtain a more accurate 
approximation of $u$ in the $L^2$ sense. 
In other words, $(P_{1,\epsilon})$ is modified as (see Section 3.1)
\begin{equation*}
  (P_{1,\varepsilon}^R):\qquad  \argmin_{\hc,~ \tensor A} \Vert \hc\Vert_1,
  \quad \text{subject to}\quad \Vert\tensor \Psi(\tensor A)\hc-\bm
  u\Vert_2\leq\epsilon, \tensor A^\trans \tensor A = \tensor I,
\end{equation*}
where $\tensor\Psi(\tensor A)$ is an $M\times N$ matrix and 
$(\tensor\Psi(\tensor A))_{ij}=\psi_j(\tensor A\bm\xi^i)$.
An alternating direction method (ADM) has been developed to iteratively identify
$\tilde{\bm c}$ and the rotation matrix $\tensor A$ based on the gradients 
of $u$. Of note, this ADM method does not guarantee to identify the
exact solution of $(P_{1,\epsilon}^R)$. It helps to identify an approximation of
$\tensor A$ such that a sparser representation of $u$ can be obtained.

In the present work, we improve the efficiency of the ADM method by using 
the sliced inverse regression (SIR) method to
provide the initial guess of the rotation matrix $\tensor A$ (Algorithm~
\ref{algo:cs_sir1}). The SIR method is used in statistics to identify important 
low-dimensional subspaces based on the training 
set $\{(\bm\xi^q, u^q)\}_{q=1}^M$. We demonstrate that the initial 
guess from SIR helps to improve the ADM algorithm accuracy in some
cases. Moreover, we propose another method that uses SIR to reduce the number 
of dimensions from $d$ to $\tilde d$, then employs ADM method to construct a 
``reduced" gPC expansion of $u$ (Algorithm~\ref{algo:cs_sir2}). In this 
case, the dimension reduction performed by 
SIR reduces the number of unknowns $N$, which can be prohibitively large for the
compressive sensing method when $d$ is large. To sum up, both new algorithms
start with SIR to identify low-dimensional subspaces. Then, this information is
used in the ADM algorithm with (Algorithm~\ref{algo:cs_sir1}) or without
(Algorithm~\ref{algo:cs_sir2}) dimension reduction to improve the accuracy of
compressive-sensing-based surrogate model construction for UQ problems. In 
this paper, we focus on problems where uncertainty (uncertain parameters) can be
described by $d$-dimensional i.i.d. Gaussian random variables 
$\bx\sim \mathcal{N}(\bm 0, \tensor I)$. This assumption is used 
broadly in physical and engineering problems, and it naturally fits the SIR 
method's requirement (see Section \ref{subsec:sir}).

The paper includes a brief review of UQ, compressive sensing methods, and the
SIR method in Section~\ref{sec:review}. Section~\ref{sec:method} describes the
proposed schemes, Algorithm~\ref{algo:cs_sir1} and Algorithm~\ref{algo:cs_sir2}.
Numerical results are presented in Section~\ref{sec:numeric}, and the 
conclusions follow in Section~\ref{sec:conclusion}.


\section{Review of compressive-sensing-based gPC and SIR methods}
\label{sec:review}
This section includes a brief review of the compressive-sensing-based gPC and 
SIR methods, which form the basis of the new method proposed 
in Section \ref{sec:method}.

\subsection{Hermite polynomial expansions}
When QoI of the problem relies on i.i.d Gaussian random variables, it can be
represented with a gPC expansion with basis functions constructed by tensor
products of univariate Hermite polynomials. Given a multi-index 
$\ba=(\alpha_1,\alpha_2,\cdots,\alpha_d), \alpha_i\in\mn\cup\{0\}$, we set
\begin{equation}\label{eq:tensor}
\psi_{\ba}(\bx) =
\psi_{\alpha_1}(\xi_1)\psi_{\alpha_2}(\xi_2)\cdots\psi_{\alpha_d}(\xi_d).
\end{equation}
A gPC expansion up to $P$-th order implies that $|\bm\alpha|\leq P$ for all
$\psi_{\bm\alpha}$ used in the expansion. For two different multi-indices 
$\ba_i=((\bm\alpha_i)_{_1}, (\ba_i)_{_2},\cdots, (\bm\alpha_i)_{_d})$ and 
$\ba_j=((\bm\alpha_j)_{_1}, (\ba_j)_{_2},\cdots, (\bm\alpha_j)_{_d})$, the Hermite
polynomials satisfy the following orthogonality condition: 
\begin{equation}
  \int_{\mr^d} \psi_{\ba_i}(\bm x)\psi_{\ba_j}(\bm x) \rho_{_{\bx}}(\bm x) \dif\bm x=
\delta_{\ba_i\ba_j} = \delta_{(\ba_i)_{_1},(\ba_j)_{_1}}
\delta_{(\ba_i)_{_2},(\ba_j)_{_2}}\cdots \delta_{(\ba_i)_{_d},(\ba_j)_{_d}},
\end{equation}
where $\delta_{(\alpha_i)_{_2},(\alpha_j)_{_2}}$ are Kronecker delta functions,
\begin{equation}
  \rho_{_{\bx}}(\bm x) =
  \rho_{_{\xi_1}}(x_1)\rho_{_{\xi_2}}(x_2)\cdots\rho_{_{\xi_d}}(x_d)
\end{equation}
and $\rho_{_{\xi_i}}(x_i)=\frac{1}{\sqrt{2\pi}}\me^{-x_i^2/2}$ because $\xi_i$
are independent Gaussian random variables. In general, when $\xi_i$ satisfy other 
distribution, $\psi_{\bm\alpha}$ can be represented as a tensor product of 
univariate polynomials associated with the PDF of $\xi_i$. In the following,
for simplicity we denote $\psi_{\ba_i}(\bx)$ as $\psi_i(\bx)$, and the gPC 
expansion used is in the form of Eq.~\eqref{eq:gpc}.


\subsection{Compressive sensing}

We first introduce the notation that denotes number of non-zeros entries in a 
vector $\bm x=(x_1,x_2,\cdots,x_N)$ \cite{Donoho06,CandesRT06,BrucksteinDE09}:
\begin{equation}
\Vert\bm x\Vert_0\Def \#\{i:x_i\neq 0\}.
\end{equation}
The vector 
$\bm x$ is called \emph{$s$-sparse} if $\Vert \bm x\Vert_0\leq s$, and $\bm x$ 
is considered a sparse vector if $s\ll N$. In practice, a very few systems have a truly 
sparse gPC coefficients $\bm c$. However, in many cases, the $\bm c$ is 
``compressible'', i.e., only a few entries make significant contribution to its
$\ell_1$ norm. Here, the $\ell_1$ norm is defined as
$\Vert\bm x\Vert_1\Def \sum_{n=1}^N |x_n|.$
Subsequently, $\bm x$ is considered sparse if $\Vert \bm x - \bm x_s\Vert_1$ is 
small for $s\ll N$, and this definition of sparsity is widely used in 
error estimation. The vector $\bm x_s$ is 
equal to $\bm x$ with all but the $s$-largest 
entries set to zero \cite{Candes08}.

The sparse vector $\bm c$ in Eq.~\eqref{eq:cs_eq} can be approximated by solving
the following $\ell_1$ minimization problem:
\begin{equation}\label{eq:lh}
(P_{1,\epsilon}): \qquad \argmin_{\hat{\bm c}}\Vert\hat{\bm c}\Vert_1, 
\text{~~subject to~~} \Vert\tensor\Psi\hat{\bm c}-\bm u\Vert_2\leq\epsilon,
\end{equation}
where $\epsilon=\Vert\bm\ve\Vert_2$. To obtain the error bound in $(P_{1,\epsilon})$, the 
\textit{restricted isometry property} (RIP) constant is introduced 
\cite{CandesT05}. For each integer $s=1,2,\cdots$, the restricted isometry constant 
$\delta_s$ of a matrix $\tensor\Psi$ is defined as the smallest number such that
\begin{equation}
(1-\delta_s)\Vert\bm x\Vert_2^2\leq \Vert\tensor\Psi\bm x\Vert_2^2\leq
(1+\delta_s)\Vert\bm x\Vert_2^2
\end{equation}
holds for all $s$-sparse vectors $\bm x$.
Cand{\`e}s et al. \cite{CandesT05} showed that if the matrix $\tensor\Psi$ satisfies 
$\delta_{2s}<\sqrt{2}-1$ (i.e., $\tensor\Psi$ satisfies ``RIP"), 
and $\Vert\bm\ve\Vert_2\leq\epsilon$, then solution 
$\hat{\bm c}$ to $(P_{1,\epsilon})$ obeys
\begin{equation}\label{eq:l1_thm}
\Vert \bm c - \hat{\bm c}\Vert_2 \leq C_1\epsilon + C_2 
\dfrac{\Vert\bm c-\bm c_s\Vert_1}{\sqrt{s}},
\end{equation}
where $C_1$ and $C_2$ are constants, and $\bm c$ is the exact vector we aim to 
approximate. This result implies that the upper bound of the error relates to the 
\emph{truncation error} and the \emph{sparsity} of $\bm c$, which is reflected in 
the first and second terms on the right-hand side of Eq. \eqref{eq:l1_thm}, respectively.

In practice, the \emph{re-weighted} $\ell_1$ minimization approach \cite{CandesWB08} 
is an improvement of the $\ell_1$ minimization method, which enhances the accuracy of 
estimating $\bm c$. It modifies $(P_{1,\epsilon})$ as
\begin{equation}\label{eq:wl1}
(P_{1,\epsilon}^W): \qquad \argmin_{\hat{\bm c}}\Vert\tensor W\hat{\bm c}\Vert_1, 
~\text{subject to}~\Vert\tensor\Psi\hat{\bm c}-\bm u\Vert_2\leq\epsilon,
\end{equation}
where $\tensor W$ is a diagonal matrix: 
$\tensor W=\text{diag}(w_1,w_2,\cdots,w_N)$. $(P_{1,\epsilon})$ can be
considered as a special case of $(P_{1,\epsilon}^W)$ with $\tensor W=\tensor I$. 
The elements $w_i$ of the diagonal matrix can be 
estimated iteratively \cite{CandesWB08, YangK13}:
in the $l$-th iteration,
 $w_i$ is set to  $w_i^{(l)}=1/(|\hat c_i^{(l-1)}|+\delta)$, where $\hat c_i^{(l-1)}$ is
the solution from the last iteration and $\hat c_i^{(0)}$ is the solution of the
standard $\ell_1$ minimization problem $(P_{1,\epsilon})$. The parameter 
$\delta>0$ is introduced to provide stability and to ensure that a zero-valued
component in $\hat{\bm c}^{(l)}$ does not prohibit a non-zero estimate at the
next step, i.e., it ensures that the weights do not become infinity.
Cand\'es et al.~\cite{CandesWB08} suggest performing two or three
iterations of this procedure. The error bound of the re-weighted $\ell_1$ 
minimization (see \cite{Needell09}) takes the same form as Eq.~\eqref{eq:l1_thm}
with different constants $C_1$ and $C_2$.
Moreover, the error term $\epsilon$ in $(P_{1,\epsilon})$ is usually not 
known \emph{a priori}, and, in the present work, we use cross-validation to 
estimate it (see the Appendix for the details). 


\subsection{Compressive-sensing-based gPC methods}
Given $M$ samples of $\bx$, we use gPC expansion Eq.~\eqref{eq:gpc} to
represent the uncertainty of QoI $u$, and we have
\begin{equation}
u(\bx^q) = \sum_{n=1}^N c_n\psi(\bx^q) + \ve(\bx^q), \quad q=1,2,\cdots,M,
\end{equation}
which can be rewritten as Eq.~\eqref{eq:cs_eq}.
A typical approach to compressive-sensing based-gPC is summarized in 
Algorithm \ref{algo:cs1} \cite{YangLBL16}.
\begin{algorithm}
\caption{Compressive-sensing-based gPC method.}
\label{algo:cs1}
\begin{algorithmic}[1]
\STATE Generate input samples $\bx^q, q=1,2,\cdots, M$ based on the distribution
of $\bx$. 
\STATE Generate output samples $u^q=u(\bx^q)$ by solving the complete model, 
e.g., running simulations, solvers, etc.
\STATE Select gPC basis functions $\{\psi_n\}_{n=1}^N$ associated with $\bx$ and
then generate the measurement matrix $\tensor\Psi$ by setting
$\Psi_{ij}=\psi_j(\bx^i)$.
\STATE Solve the optimization problem $(P_{1,\epsilon})$:
\[\argmin_{\hat{\bm c}}\Vert\hat{\bm c}\Vert_1, ~
\text{subject to} \Vert\tensor\Psi\hat{\bm c}-\bm u\Vert_2\leq\epsilon,\]
where $\bm u=(u^1,u^2,\cdots,u^M)^\trans$, and $\epsilon$ is obtained by
cross-validation (see Algorithm~\ref{algo:cross} in Appendix). If 
the re-weighted $\ell_1$ method is employed, solve $(P_{1,\epsilon}^W)$ instead.
\STATE Construct the gPC expansion as 
$u(\bx)\approx \sum_{n=1}^N \hat c_n\psi_n(\bx)$.
\end{algorithmic}
\end{algorithm}

\subsection{Sliced Inverse Regression}
\label{subsec:sir}
SIR is an effective approach for seeking the important subspaces in the parameter
space \cite{Li91}. As an illustration, consider 
$u(\bm\xi)=u(\xi_1,\xi_2) = (\xi_1+\xi_2)^2$. Then, if we define 
$\hat{\tensor A}=(1/\sqrt{2}, 1/\sqrt{2})$ and
$\eta_1=\hat{\tensor A}\bm\xi$, $u$ only depends on $\eta_1$.
Consequently, the dimension is reduced from $d=2$ to reduced 
dimension $\tilde d=1$ because we only need one input random variable to fully
capture the statistical property of $u$.
Unlike the outer product gradients (OPGs) \cite{HardinH07, Xia07} or active subspace method
\cite{ConstantineDW14} where gradients information is used to identify
$\hat{\tensor A}$, the SIR method uses conditional expectation $\mexp{\bm\xi | u}$.
$\mexp{\bm\xi | u}$ is a $d$-dimensional random vector because $u$ is
random. As $u$ varies, $\mexp{\bm\xi | u}$ draws a curve in the parameter space, 
which is called \emph{inverse regression curve}. It has been shown that this curve 
resides in the desired subspace for dimension reduction (named \emph{central subspace}) 
if $\bm\xi$ follows an elliptically symmetric distribution \cite{Li91}, e.g., 
the multivariate Gaussian distribution. Based on this property, we choose
the matrix $\hat{\tensor A}$ such that its columns consist of the eigenvectors 
corresponding to the non-zero eigenvalues of the covariance matrix 
$\tensor V=\text{var}\left\{\mexp{\bm\xi | u}\right\}$. An estimate of $\hat{\tensor A}$
is summarized in Algorithm \ref{algo:sir}, originally proposed in \cite{Li91}. A software 
package implementing the algorithm is available in \cite{Weisberg02}.
\begin{algorithm}
  \caption{Sliced inverse regression algorithm.}
\label{algo:sir}
\begin{algorithmic}[1]


\STATE Generate i.i.d. samples of input parameters $\bm\xi^i$, $i=1,...,M$, and
compute the corresponding values of QoI $u^i$.

\STATE Divide the range of $u^i$, i.e., $[\min(u^i),\max(u^i)]$, into $H$
  non-overlapping slices, $J_1,...,J_H$: $[u_0, u_1), [u_1, u_2), \cdots, [u_{H-1}, u_H]$,
    where  $\min(u^i)=u_0<u_1<u_2<\cdots<u_{H-1}<u_H=\max(u^i)$,
each containing approximately an equal number of data points.

\STATE Compute the within-slice mean of $\bm\xi$ over each slice, which is a
crude estimate of the conditional expectation $\mexp{\bm\xi| u}$:

\[ \bar{\bm\xi}_h = \frac{1}{n_h} \sum_{u^i \in J_h} \bm\xi^i, \; h=1,...,H, \]

where $n_h$ is the number of data points falling in the $h$th slice.

\STATE Compute the $d \times d$ matrix 
\[ \tensor V =  \sum_{i=1}^H \frac{n_h}{n} \bar{\boldsymbol{\xi}}_h \bar{\boldsymbol{\xi}}_h^\trans. \]
This is the sample estimate of the covariance matrix of the random vector $E (
\boldsymbol\xi | \tilde u)$, where $\tilde u = \sum_{h=1}^H h I ( u \in J_h)$
and $I(\cdot)$ is the indicator function.


\STATE Compute the eigen-decomposition of $\tensor V$: $\tensor V=\tensor
U_V\tensor \Lambda_V\tensor U_V^\trans$, where $\tensor\Lambda_V$
is a diagonal matrix consists of eigenvalues: $(\tensor\Lambda_{V})_{ii}=(\lambda_V)_i$ 
($1\leq i\leq d$) with $(\lambda_V)_1\geq(\lambda_V)_2\geq\cdots \geq (\lambda_V)_d\geq 0$ 
and $\tensor U_V\tensor U_V^\trans=\tensor I$.

\STATE Report the estimated transformation matrix as a submatrix
consisting of the first $\tilde d$ (reduced dimension) rows of $\tensor U_V^\trans$:
$\hat{\tensor A}=[(\tensor U_V)_1, ..., (\tensor U_V)_{\tilde d}]^\trans$, where $\tilde d\leq d$.

\end{algorithmic}
\end{algorithm}
Of note, most applications of SIR are concerned with dimension reduction by
choosing $\tilde d$ to be as small as possible, i.e., smaller than $d$  
(see Algorithm \ref{algo:sir}). In this work, we use this setting in 
Algorithm \ref{algo:cs_sir2}. On the other hand, we use $\hat{\tensor A}$ with 
$\tilde d=d$ to obtain an initial guess for the ADM algorithm (see 
Algorithm \ref{algo:cs_sir1}). 

Notably, SIR can be considered as an approach within the framework of \emph{sufficient
dimension reduction} (SDR). To simplify the model $u(\bm\xi)$, an effective 
modeling strategy is to assume that only a few subspaces make major 
contributions to $u$. A formal definition tailored from \cite{Li91} 
in \cite{LiLinLi16} is as follows:

\noindent\textbf{Definition:} Given the $d$-dimensional model $u(\bm\xi)$,
a \textit{dimension reduction} is a mapping from the $d$-dimensional input to a
$\tilde d$-dimensional vector, $\boldsymbol{\eta} = \tensor A \boldsymbol{\xi}$, 
where $\tensor A \in \mr^{\tilde d \times d}$, $\tilde d < d$, $\tensor A
\tensor A^\trans = \tensor I$ is the identity matrix. A dimension reduction is 
\textit{sufficient} if the following equation holds for any $\boldsymbol{\xi} \in \mr^d$:
\[u(\boldsymbol{\xi}) = u(\tensor A^\trans \tensor A \boldsymbol{\xi} )
\equiv f(\tensor A^\trans \boldsymbol{\eta}).\]
In other words, $u$ only relies on $\tilde d$ variables $\eta_1, \eta_2,
\cdots, \eta_{\tilde d}$. For example, the active subspace method
\cite{ConstantineDW14}, basis adaptation method \cite{TipG14}, and SIR method
aim to identify this low-dimensional structure by computing $\tensor A$ in a
different manner.


\section{SIR-aided Rotated Compressive Sensing Method}
\label{sec:method}

This section details two new applications of the SIR method for the 
compressive-sensing-based gPC, which is the main contribution of this work. 

\subsection{Alternating direction method for increasing sparsity}
\label{subsec:method}

In practical problems, if the truncation error $\epsilon$ is sufficiently small,
then the second term on the right-hand side of Eq.~\eqref{eq:l1_thm} dominates
the upper bound of the error. Hence, to improve the accuracy of the gPC
expansion, we need to decrease 
$\Vert\bm c-\bm c_s\Vert_1/\sqrt{s}$. Our goal is to seek $\tensor A$ (and
$\bm\eta=\tensor A\bm\xi$) such that in Eq.~\eqref{eq:gpc2},
$\Vert\tilde{\bm c}-\tilde{\bm c}_s\Vert_1<\Vert\bm c-\bm c_s\Vert_1$. 
In other words, we rewrite the standard $\ell_1$ minimization problem
$(P_{1,\epsilon})$ as 
\begin{equation}\label{eq:rot_l1}
  (P_{1,\varepsilon}^R): \qquad \argmin_{\hc,~ \tensor A} \Vert \hc\Vert_1,
  \quad \text{subject to}\quad \Vert\tensor \Psi(\tensor A)\hc-\bm u\Vert_2\leq\epsilon,
  \tensor A^\trans \tensor A = \tensor I.
\end{equation}
where $\tensor\Psi(\tensor A)$ is a matrix and 
$(\tensor\Psi(\tensor A))_{ij}=\psi_j(\tensor A\bm\xi^i)$. In the ADM algorithm
proposed in \cite{YangLBL16}, gradient information is used to identify 
$\tensor A$. A ``gradient matrix" is defined as
\begin{equation}\label{eq:grad_mat}
\tensor G \Def
\mexp{\nabla u(\bx)\otimes\nabla u(\bx)} = \tensor U\tensor \Lambda\tensor
U^\trans,
  \quad \tensor U\tensor U^\trans = \tensor I,
\end{equation}
where $\tensor G$ is symmetric, 
$\nabla u(\bx) = (\partial u/\partial \xi_1,
\cdots, \partial u/\partial \xi_d)^\trans$ is a column vector, 
$\tensor U=(\bm U_1,\cdots, \bm U_d)$ is an orthogonal matrix
consisting of eigenvectors $\bm U_i$, and 
$\tensor\Lambda=\text{diag}(\lambda_1,\cdots,\lambda_d)$ with 
$\lambda_1\geq\lambda_2\geq\cdots\lambda_d\geq 0$ is a diagonal matrix with elements 
representing variation of the system along the respective 
eigenvectors. Then, $\tensor A$ can be chosen as the unitary matrix $\tensor
U^\trans$, which defines a rotation in $\mr^d$  projecting $\bx$ on the 
eigenvectors $\bm U_i$. If only a few $\lambda_i$s are very large (compared 
with other $\lambda_i$s), the rotation that maps $\bm\xi$ to 
$\bm\eta=\tensor A\bm\xi$ helps to concentrate the dependence of $u$ primarily
on those few new random variables $\eta_i$ due to the larger variation of $u$ 
along the directions of the corresponding eigenvectors. Therefore, the resulting
coefficients $\tilde{\bm c}$ can be sparser than $\bm c$. This approach of 
constructing $\tensor G$ from active subspace (proposed in \cite{ConstantineDW14})
is similar to the method of OPGs in statistics \cite{HardinH07, Xia07}. The 
gradient of $u$ also has been used to improve the efficiency of compressive 
sensing in the \emph{gradient-enhanced} method \cite{LiARH11,JakemanES14,PengHD15}.

Because the explicit form of $u$ or $\nabla u$ is unknown, an ADM algorithm is
proposed to identify $\tensor A$ and $\tilde c$ iteratively. As noted in the
introduction, this work aims to solve $(P_{1,\epsilon}^R)$ when $\xi_i$ are i.i.d.
Gaussian random variables. Therefore, we use the algorithm from
\cite{YangLBL16} that is summarized in Algorithm \ref{algo:cs_rot}. A general
form of this algorithm that handles $\bm\xi$ of different distributions can be
found in \cite{YangWL17}, but it is beyond the scope of this work.
\begin{algorithm}[h]
  \caption{Alternating direction method for solving $(P_{1,\epsilon}^R)$ when
  $\xi_i$ are i.i.d. Gaussian random variables.}
\label{algo:cs_rot}
\begin{algorithmic}[1]
\STATE Generate input samples $\bx^q, q=1,2,\cdots, M$ based on the distribution
of $\bx$. 
\STATE Generate output samples $u^q=u(\bx^q)$ by solving the complete model, 
e.g., running simulations, solvers, etc.
\STATE Select gPC basis functions $\{\psi_n\}_{n=1}^N$ as normalized Hermite polynomials and
then generate the measurement matrix $\tensor\Psi$ by setting
$\Psi_{ij}=\psi_j(\bx^i)$.
\STATE Solve the optimization problem $(P_{1,\epsilon})$:
\[\arg \min_{\hat{\bm c}}\Vert\hat{\bm c}\Vert_1, ~
\text{subject to}~ \Vert\tensor\Psi\hat{\bm c}-\bm u\Vert_2\leq\epsilon.\]
\STATE Set counter $l=0$, $\eta^{(0)}=\bx$, $\tilde{\bm c}^{(0)}=\hat{\bm c}$,
compute $\tensor K_{ij}, i,j=1,2,\cdots,d$.
\STATE $l=l+1$. Construct ${\tensor G}^{(l)}$ as 
$G^{(l)}_{ij}=(\tilde{\bm c}^{(l-1)})^T \tensor K_{ij} \tilde{\bm c}^{(l-1)},
i,j=1,2,\cdots, d$.
Then, compute eigen-decomposition of $\tensor G^{(l)}$:
\[\tensor G^{(l)}= 
\tensor U^{(l)}\tensor\Lambda^{(l)}\left(\tensor U^{(l)}\right)^T.\]
\STATE Set $\bm\eta^{(l)}=\left(\tensor U^{(l)}\right)^\trans\bm\eta^{(l-1)}$.
Then compute samples $(\bm\eta^{(l)})^q=\left(\tensor U^{(l)}\right)^\trans(\bm\eta^{(l-1)})^q, q=1,2,\cdots,M$.
Also, construct the new measurement matrix $\tensor\Psi^{(l)}$ with
$\Psi^{(l)}_{ij}=\psi_j\left((\bm\eta^{(l)})^i\right)$. 
\STATE Solve the optimization problem $(P_{1,\epsilon^{(l)}})$:
\[\arg \min_{\hat{\bm c}}\Vert\hat{\bm c}\Vert_1, \quad\text{subject to}~
\Vert\tensor\Psi^{(l)}\hat{\bm c}-\bm u\Vert_2\leq\epsilon^{(l)},\]
and set $\tilde{\bm c}^{(l)}=\hat{\bm c}$. 
\STATE If $|\Vert\tensor U^{(l)}\Vert_1-d|<\theta$, where the 
threshold $\theta$ is a positive real number, then stop the iterations. 
Otherwise, go to Step 6.
\STATE Set 
\[\tensor A^{(l)}=\left(\tensor U^{(1)}\tensor U^{(2)}\cdots\tensor U^{(l)}\right)^T\]
and construct gPC expansion as $u(\bx)\approx u_g(\bm\xi)=v_g(\bm\eta^{(l)})=\sum_{n=1}^N \tilde{c}^{(l)}_n\psi_n(\tensor A^{(l)}\bx)$.
\end{algorithmic}
\end{algorithm}
The matrix $\tensor K_{ij}$ in Step 5 is defined as
\begin{equation}\label{eq:kernel}
(\tensor K_{ij})_{kl} = \mexp{\dfrac{\partial\psi_k(\bm\xi)}{\partial\xi_i}
                \cdot\dfrac{\partial\psi_{l}(\bm\xi)}{\partial\xi_j}}, \quad
                1\leq k,l\leq N.
\end{equation}                
The analytic form of $\tensor K_{ij}$ is 
\begin{equation}
(\tensor K_{ij})_{kl}  =
 \sqrt{(\bm\alpha_k)_{_i}(\bm\alpha_l)_{_j}}
\delta_{(\bm\alpha_k)_{_i}-1,(\bm\alpha_l)_{_i}}
\delta_{(\bm\alpha_k)_{_j},(\bm\alpha_l)_{_j}-1}\cdot
\prod_{\substack{m=1\\ m\neq i,m\neq j}}\delta_{(\bm\alpha_k)_{_m},(\bm\alpha_l)_{_m}}.
\end{equation}
Algorithm~\ref{algo:cs_rot} takes advantage of the Gaussian random
variables properties in the following ways:
in each iteration, $\bm\eta$ is updated as 
$\bm\eta^{(l)}=\left(\tensor U^{(l)}\right)^\trans\bm\eta^{(l-1)}$ in Step 7,
and both $\bm\eta^{(l)}$ and $\bm\eta^{(l-1)}$ follow the Gaussian
distribution $\mathcal{N}(\bm 0, \tensor I)$ because it is a orthogonal matrix.
Therefore, we only need a ``correction" of $\tensor A$ (i.e., $\tensor U^{(l)}$) 
in each iteration, and the matrix $\tensor A$ is computed after all iterations 
are completed in Step 10. More specifically, in each iteration, $\tensor A^{(l)}$ 
can be computed as $\tensor A^{(l)}=\left(\tensor U^{(l)}\right)^\trans \tensor A^{(l-1)}$,
but $\tensor A^{(l)}$ is not needed explicitly to update $\bm\xi^{(l)}$. 
Moreover, in Step 8, $\epsilon^{(l)}$ may vary in different iterations. It is
usually sufficient to test two or three different values on the interval 
$[\epsilon/5,\epsilon]$ using cross-validation to identify $\epsilon^{(l)}$. 
The stopping criterion in Step 9 measures the distance between $\tensor U^{(l)}$ 
and the identity or permutation matrix \cite{YangLBL16}.
Empirically, the threshold $\theta$ can be taken as 
$0.2d\sim 0.3d$ when the dimension $d$ is $\mo(10)$ and $0.5d\sim 0.8d$ when $d$
is $\mo(100)$.

\subsection{SIR-aided ADM for increasing sparsity}
\label{subsec:sir1}
The first proposd approach involves using SIR to provide an initial guess for the
aforementioned ADM algorithm, and improve an estimate of the rotational matrix 
$\tensor A$. Specifically, in Algorithm \ref{algo:cs_rot}, the iteration 
starts with initial guess $\tilde{\bm c}^{(0)}$ obtained at Step 4. Then, the 
initial guess of $\tensor A$ is constructed based on $\tilde{\bm c}^{0}$ in 
Step 6. Instead, we can start with an initial guess of $\tensor A$ from SIR and
compute $\tilde c^{(1)}$. In this approach, we do not solve
$(P_{1,\epsilon})$ to provide an initial guess of $\tilde{\bm c}$, i.e., we skip
Step 4\ in Algorithm~\ref{algo:cs_rot}. The new algorithm--\emph{SIR-based
ADM}(SADM)--is summarized in Algorithm~\ref{algo:cs_sir1}.
\begin{algorithm}[h]
\caption{Alternating direction method of solving $(P_{1,\epsilon}^R)$ based on
SIR (SADM) when $\xi_i$ are i.i.d. Gaussian random variables.}
\label{algo:cs_sir1}
\begin{algorithmic}[1]
\STATE Generate input samples $\bx^q, q=1,2,\cdots, M$ based on the distribution
of $\bx$. 
\STATE Generate output samples $u^q=u(\bx^q)$ by solving the deterministic
problem with input $\bx^q$.
\STATE Select gPC basis functions $\{\psi_n\}_{n=1}^N$ as normalized Hermite polynomials.
\STATE Run Algorithm~\ref{algo:sir} with the training set 
$\{(\bm\xi^q, u^q)\}_{q=1}^M$, to obtain $\hat{\tensor A}$ by setting $\tilde
d=d$, then set $\tensor U^{(1)}=\hat{\tensor A}^\trans$.
\STATE Set $\eta^{(0)}=\bx$ and counter $l=1$. Then compute $\tensor K_{ij},
i,j=1,2,\cdots,N$.
\STATE Set $\bm\eta^{(l)}=\left(\tensor U^{(l)}\right)^\trans\bm\eta^{(l-1)}$.
Then, compute samples $(\bm\eta^{(l)})^q=\left(\tensor U^{(l)}\right)^\trans(\bm\eta^{(l-1)})^q, q=1,2,\cdots,M$.
Also, construct the measurement matrix $\tensor\Psi^{(l)}$ as
$\Psi^{(l)}_{ij}=\psi_j\left((\bm\eta^{(l)})^i\right)$. 
\STATE Solve the optimization problem $(P_{1,\epsilon^{(l)}})$:
\[\arg \min_{\hat{\bm c}}\Vert\hat{\bm c}\Vert_1, \quad\text{subject to}~
\Vert\tensor\Psi^{(l)}\hat{\bm c}-\bm u\Vert_2\leq\epsilon^{(l)},\]
and set $\tilde{\bm c}^{(l)}=\hat{\bm c}$. 
\STATE If $|\Vert\tensor U^{(l)}\Vert_1-d|<\theta$, where the 
threshold $\theta$ is a positive real number, then stop the iterations. Otherwise, 
Set $l=l+1$ and construct ${\tensor G}^{(l)}$ as 
$G^{(l)}_{ij}=(\tilde{\bm c}^{(l-1)})^\trans \tensor K_{ij} \tilde{\bm c}^{(l-1)},
i,j=1,2,\cdots, d$.
Then compute eigen-decomposition of $\tensor G^{(l)}$:
\[\tensor G^{(l)}= 
\tensor U^{(l)}\tensor\Lambda^{(l)}\left(\tensor U^{(l)}\right)^\trans,\]
and go to Step 6.
\STATE Set 
\[\tensor A^{(l)}=\left(\tensor U^{(1)}\tensor U^{(2)}\cdots\tensor U^{(l)}\right)^T,\]
and construct gPC expansion as 
$u(\bx)\approx u_g(\bm\xi)=v_g(\bm\eta^{(l)})=\sum_{n=1}^N \tilde{c}^{(l)}_n\psi_n(\tensor A^{(l)}\bx)$.
\end{algorithmic}
\end{algorithm}

The  difference between Algorithms~\ref{algo:cs_rot} and
\ref{algo:cs_sir1} is the initial guess of $\tensor A$, i.e., how to compute 
$\tensor U^{(1)}$. Section~\ref{sec:numeric} shows how the initial 
guess provided by SIR yields a more accurate estimate of $u_g$ in
our test cases. In the compressive sensing theory, there is a requirement on 
the size $M$ of available data for high probability of 
the signal recovery. For example, an $s$-sparse (univariate) trigonometric
polynomial of maximal degree $P$ (i.e., $N=P+1$) can be recovered from
$M\asymp s\log^4(P)$ sampling points \cite{CandesRT06, RudelsonV08}, and an 
$s$-sparse (univariate) Legendre polynomial of maximal degree $P$ 
(again, $N=P+1$) can be recovered from $M\asymp s\log^3(s)\log(P)$ sampling
points from a Chebyshev measure \cite{RauhutW12}. 
If the number of sampling points is too small (compared with $N$), 
there is no guarantee that the compressive sensing results will be accurate
even if $s$ is small.
For example, we limit the sample size $M$ as $\mathcal{O}(100)$
in the numerical tests (Section \ref{sec:numeric}), which is
typical in practical problems. When dimension $d$ is high, $N$ becomes very 
large, thus, the standard compressive sensing method may not work well.
In this scenario, the gradient computed from a truncated gPC expansion
may not provide the optimal initial guess for the ADM algorithm 
because $\tilde{\bm c}$ is inaccurate. Unlike the compressive sensing method,
which is based on a regression form of $u$, the
SIR method does not assume a specific regression form of $u$,
nor does it use gradient
information to identify $\tensor A$. A theoretical analysis in \cite{LiLinLi16}
demonstrates that if there is an ``optimal" rotation matrix $\tensor A$
(e.g., this matrix exisits in numerical example 4.1), then
$\Vert\hat{\tensor A}-\tensor A\Vert_2$ is $\mo(M^{-1})$, where 
$\hat{\tensor A}$ is found from SIR. Notably, $N$ is not explicitly included 
in this estimate because SIR does not assume an expansion form of $u$. 

\subsection{SIR-aided alternating direction method based on dimension reduction}
\label{subsec:sir2}
The second proposed approach is to precede the ADM algorithm with dimension 
reduction by SIR. 
As noted in the discussion regarding sample size requirement, when $M$ is much
smaller than $N$, even Algorithm \ref{algo:cs_sir1} may not be directly 
applicable because the compressive sensing algorithm cannot provide accurate
results. As an alternative, we propose a new algorithm combining compressive
sensing and dimension reduction performed by SIR (SADMDR). The idea of this
method is to reduce $N$ before using the ADM algorithm.

In a PC expansion of $u$ up to a polynomial order $P>1$, $N$ grows exponentially
with increasing $d$. Consequently, in problems with large $d$, $M$ could be much
smaller than $N$, and the compressive sensing results are expected to be less  
accurate. Although, a larger $d$ also effects the accuracy of the result from
SIR, SIR is still expected to provide a better initial guess than compressive
sensing in this scenario. 
Furthermore, if we keep all Hermite polynomials up to order $P$ in the gPC
expansion, the number of unknown coefficients in $u_g(\tilde\bx)$ is 
$N=\bigl(\begin{smallmatrix} P+ d \\ P \end{smallmatrix}\bigr)$. Given the 
limited available data $\{(\bm\xi^q,u^q)\}_{q=1}^M$, $N$ cannot be too large, 
otherwise the compressive sensing method would not provide an accurate estimate 
of the gPC expansion. This implies that $P$ should be small (in many cases no
larger than $2$) if $d$ is large and $M$ is small. However, gPC expansions
with small $P$ could be inaccurate, especially for estimating second-and 
higher-order moments of $u$. For example, the variance of $u$ is estimated as 
$\text{Var}\left\{u\right\}\approx
\text{Var}\left\{u_g\right\}=\sum_{n=2}^Nc_n^2$.
An accurate estimate of $\text{Var}\left\{u\right\}$ requires a sufficient number
of higher-order Hermite polynomials in the expansion of $u_g$. To some extent, 
the SIR method can help to solve this dilemma. Originally, SIR was designed for 
dimension reduction, with $\tilde d$ chosen smaller (in many cases, much smaller)
than $d$ in Step 7 of Algorithm~\ref{algo:sir}. In other words, the $d$-dimensional
vector $\bm\xi$ is projected to a $\tilde d$-dimensional vector 
$\tilde{\bm\xi}=\hat{\tensor A}\bm\xi$, where $\tilde d<d$. The reduced vector
$\tilde{\bm\xi}$ can be used to construct a ``reduced'' gPC expansion 
$\tilde u_g$ as an approximation of $u$, such that $\tilde u_g$ has
approximately the same statistical properties (mean, standard deviation, PDF,
etc) as $u$. Because $d$ is reduced to $\tilde d$, it is possible to use larger
$P$ in the gPC expansion while keeping $N$ in an appropriate range,
allowing the compressive sensing method to obtain an accurate approximation of
$u$ with $M$ sampling points. The SADMDR method is described in 
Algorithm~\ref{algo:cs_sir2}. 
\begin{algorithm}[h]
  \caption{Alternating direction method based on dimension reduction by SIR
  (SADMDR) when
$\xi_i$ are i.i.d. Gaussian random variables.}
\label{algo:cs_sir2}
\begin{algorithmic}[1]
\STATE Generate input samples $\bx^q, q=1,2,\cdots, M$ based on the distribution
of $\bx$. 
\STATE Generate output samples $u^q=u(\bx^q)$ by solving the complete model, 
e.g., running simulations, solvers, etc.
\STATE Run Algorithm~\ref{algo:sir} with the training set 
$\{(\bm\xi^q, u^q)\}_{q=1}^M$, to obtain $\hat{\tensor A}$ by setting 
$\tilde d<d$. Then set $\tilde{\bm\xi}=\hat{\tensor A}\bm\xi$ and compute
corresponding sampling points $\tilde{\bm\xi}^q=\hat{\tensor A}\bm\xi^q,
q=1,2,\cdots,M$.
\STATE Run Algorithm~\ref{algo:cs_rot} based on training sets
$\{(\tilde{\bm\xi}^q, u^q)\}_{q=1}^M$ to obtain $\tilde u_g$ as
\[\tilde u_g(\tilde{\bm\xi})=\tilde v_g(\tilde{\bm\eta}^{(l)})=\sum_{n=1}^N
\tilde{\tilde{c}}^{(l)}_n\psi_n(\tilde{\tensor A}^{(l)}\tilde\bx).\]
\end{algorithmic}
\end{algorithm}

In the $\tilde u_g$ expansion, $N$ is set as
$\bigl(\begin{smallmatrix} P+ \tilde d \\ P \end{smallmatrix}\bigr)$, and it is
possible to use larger $P$ than in the original problem because $\tilde d$ is 
smaller than $d$. 
By reducing dimensionality, we lose some information about $u$, unless this is 
a sufficient dimension reduction, i.e., the system has a lower dimensional
representation (see the remark at the end of Section
\ref{sec:review}). An example of SDR,  
$u(\bm\xi)=u(\xi_1,\xi_2) = (\xi_1+\xi_2)^2$, is shown in Section 2.4. Here, $u$
only depends on $\xi_1+\xi_2$. Therefore, setting 
$\tilde d=1, \hat{\tensor A}=(1/\sqrt{2},1/\sqrt{2})$ does not lead to any loss 
of information about $u$. However, this is not true for most practical problems, 
and truncating the dimension too aggressively, no matter how large $P$ is used in 
the ${\tilde u}_g$ expansion, usually leads to a poor approximation of $u$. On
the other hand, choosing a relatively large $\tilde d$ (and possibly keeping 
``unimportant'' information about $u$) requires selecting $P$ that is too large
for the compressive sensing method to produce an accurate estimate of $u$. In 
practice, the value of $\tilde d$ is decided based on the change of magnitude of
the eigenvalues $(\lambda_V)_i$ in SIR. For example, one can select $\tilde d$ 
such that $\sum_{i=1}^{\tilde d} (\lambda_V)_i \geq a\sum_{i=1}^d (\lambda_V)_i$
and $a<1$. In this work, we use an R package implementation of
SIR with a $p$-value test to determine $\tilde d$ \cite{Weisberg02}. After
$\tilde d$ is found, we select $P$ such that $N$ is between $2M$ to $5M$. This
choice of $P$ usually ensures that the compressive sensing method will produce
an accurate approximation of $u$. If a \emph{prior} knowledge of the sparsity 
of $u$ is available, $P$ can be selected more appropriately. This may require 
a numerical analysis of the partial differential equation (PDE), domain 
knowledge of the system, etc., and it is beyond the scope of our work.

\noindent\textbf{Remark:} In this study, we focus on Gaussian random variables
and Hermite polynomials. Algorithm~\ref{algo:cs_rot} can also be applied to other type of
random variables and their associated orthogonal polynomials \cite{YangWL17}. 
The SIR works well for Gaussian random variables but performs much worse when
the distribution of random variables deviates from the Gaussian case. There are
several methods designed for more general cases, e.g. sliced average variance
estimator \cite{CookW91}, minimum average variance estimator \cite{Xia02}.
This is beyond the discussion of this work and we refer the interested readers 
to these literatures.


\section{Numerical Examples}
\label{sec:numeric}

In this section, five numerical examples are used to demonstrate the effectiveness
of the proposed method. In examples 1 and 2, the test functions are $12$- and 
$20$-dimensional polynomials, respectively, and Algorithm \ref{algo:cs_sir1} is
used to construct the Hermite polynomial expansion $u_g$. The accuracy of different 
methods is measured by the relative $L_2$ error: $(\Vert u - u_g\Vert_2)/\Vert u\Vert_2$. 
The integral in
\begin{equation}
\Vert u(\bx)\Vert_2 = \left(\int_{\mr^d} u(\bm x)^2\rho_{\bm\xi}(\bm x)\dif\bm x\right)^{1/2}
\end{equation}
and $\Vert u-u_g\Vert_2$ is approximated with a high-level sparse grids method 
based on one-dimensional Gauss quadrature and the Smolyak structure
\cite{Smolyak63} to guarantee accurate numerical integration.  
Examples 3-5 are high-dimensional (from $100$- to $500$-dimensional)
stochastic PDE problems and a polynomial test function,
and the relative errors of the mean and standard deviation are presented to 
compare the accuracy of different methods. In these three examples, the reference
solution of the mean and standard deviation of $u$ are obtained from $10^6$ Monte 
Carlo (MC) realizations. All relative errors presented in this section are obtained
from $100$ independent replicates for each sample size $M$. Namely, we generate $100$
independent sets of input samples $\bx^q,q=1,2,\cdots,M$, compute corresponding
relative errors, and report the average of these error samples using symbols. In
the first two examples, we investigate the relative error of various methods as a
function of the available data size relative to the number of unknowns, e.g., 
$M/N$. In examples 3-5, we compute error as a function of $M$, and we also 
present the quantiles (25th and 75th percentiles) using horizontal bars. We use
the MATLAB package \texttt{SPGL1} \cite{BergF08,spgl1} to solve $(P_{1,\epsilon})$
and the R package \texttt{dr} implementation of SIR \cite{Weisberg02}. 

\subsection{Ridge function}
\label{subsec:ex1}
\begin{figure}[h]
\centering
\includegraphics[width=0.45\textwidth]{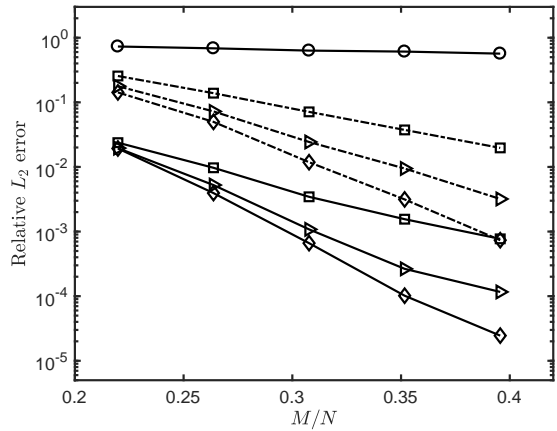}\quad
\caption{Results for the ridge function. ``$\circ$": standard $\ell_1$,
``$\square$": $\ell_1$ with 3 rotations, ``$\triangleright$": $\ell_1$ with 6 
rotations, ``$\diamond$": $\ell_1$ with 9 rotations. Dash lines result from 
using initial guess of $\tensor A$ based on compressive sensing results
(Algorithm~\ref{algo:cs_rot}), while
solid lines are results using an initial guess of $\tensor A$ from SIR
(Algorithm~\ref{algo:cs_sir1})}.
\label{fig:ex1_rmse}
\end{figure}

Consider the following ridge function:
\begin{equation}\label{eq:ex1}
u(\bx) = \sum_{i=1}^d \xi_i + 0.25\left(\sum_{i=1}^d \xi_i\right)^2
       + 0.025\left(\sum_{i=1}^d \xi_i\right)^3.
\end{equation}
This example is used in \cite{YangLBL16, YangWL17} to demonstrate the 
effectiveness of the iterative rotational $\ell_1$ compressive sensing method. 
This ridge function is unique in that all $\xi_i$ are equally important. 
Hence, adaptive methods that build the surrogate model hierarchically based on 
the importance of each $\xi_i$ (e.g., \cite{YangCLK12, ZhangYMMEKD14}) may not 
be efficient--or even work at all. The Hermite polynomial expansion of
$u(\bm\xi)$ with $P=3$ is not exactly sparse as none of the coefficients are
zero. The rotation matrix 
\begin{equation}\label{eq:ex_rot}
\tensor A = 
\begin{pmatrix}
d^{-1/2} & d^{-1/2} & \cdots & d^{-1/2} \\
 &   &  & \\
 & \multicolumn{2}{c}{\tilde{\tensor{A}}}  & \\
 &   &  & 
\end{pmatrix},
\end{equation}
reduces $u$ to a concise form:
\[u(\bx)=u(\bm\eta)=d^{1/2}\eta_1+0.25d\eta_1^2+0.025d^{3/2}\eta_1^3,\]
where $\tilde{\tensor A}$ is a $(d-1)\times d$ matrix chosen to ensure that
$\tensor A$ is orthonormal and  
$\eta_1=(\sum_{i=1}^d\xi_i)/d^{1/2}$.
If the set of the basis functions remains unchanged, all of the
polynomials not related to $\eta_1$ make no contribution to the expansion of $u$,
which implies that we obtain an $s$-sparse Hermite polynomial expansion with
$s=4$. Specifically, only four Hermite polynomials (from zero-th
order term to the third-order term) are needed to represent $u(\bm\eta)$. We
demonstrated in \cite{YangLBL16, YangWL17} that the alternating direction method
is able to detect the optimal structure using iterations and yield an accurate 
approximation of $u$. Here, we repeat the same test by setting $d=12$
(hence, $N=455$ for $P=3$) to demonstrate the effectiveness of the new method 
and to compare it with the result in \cite{YangLBL16}. Figure~\ref{fig:ex1_rmse}
represents the relative errors. Clearly, the standard $\ell_1$ 
minimization is not effective as the relative error is more than $50\%$ even 
when $M/N$ is close to $0.4$. Also, it is demonstrated in \cite{YangLBL16} that the
re-weighted $\ell_1$ does not help in this case. The ADM
Algorithm~\ref{algo:cs_rot} (dash lines) improves the accuracy by up to two 
magnitudes using $9$ iterations, and the new Algorithm~\ref{algo:cs_sir1} (solid
lines) is able to further improve the accuracy by one more magnitude. Comparing
results denoted by the same symbols (triangles, squares, and diamonds) on the 
dash lines (Algorithm~\ref{algo:cs_rot}) and  solid lines 
(Algorithm~\ref{algo:cs_sir1}) shows that for each fixed $M/N$, the symbols
on the solid lines are one magnitude lower than those on the corresponding dash
lines. These results demonstrate that using an initial guess of $\tensor A$ from
SIR improves the accuracy of the $u_g$ approximation of $u$.


\subsection{Function with high compressibility}
\label{subsec:ex2}
\begin{figure}[h]
\centering
\includegraphics[width=0.46\textwidth]{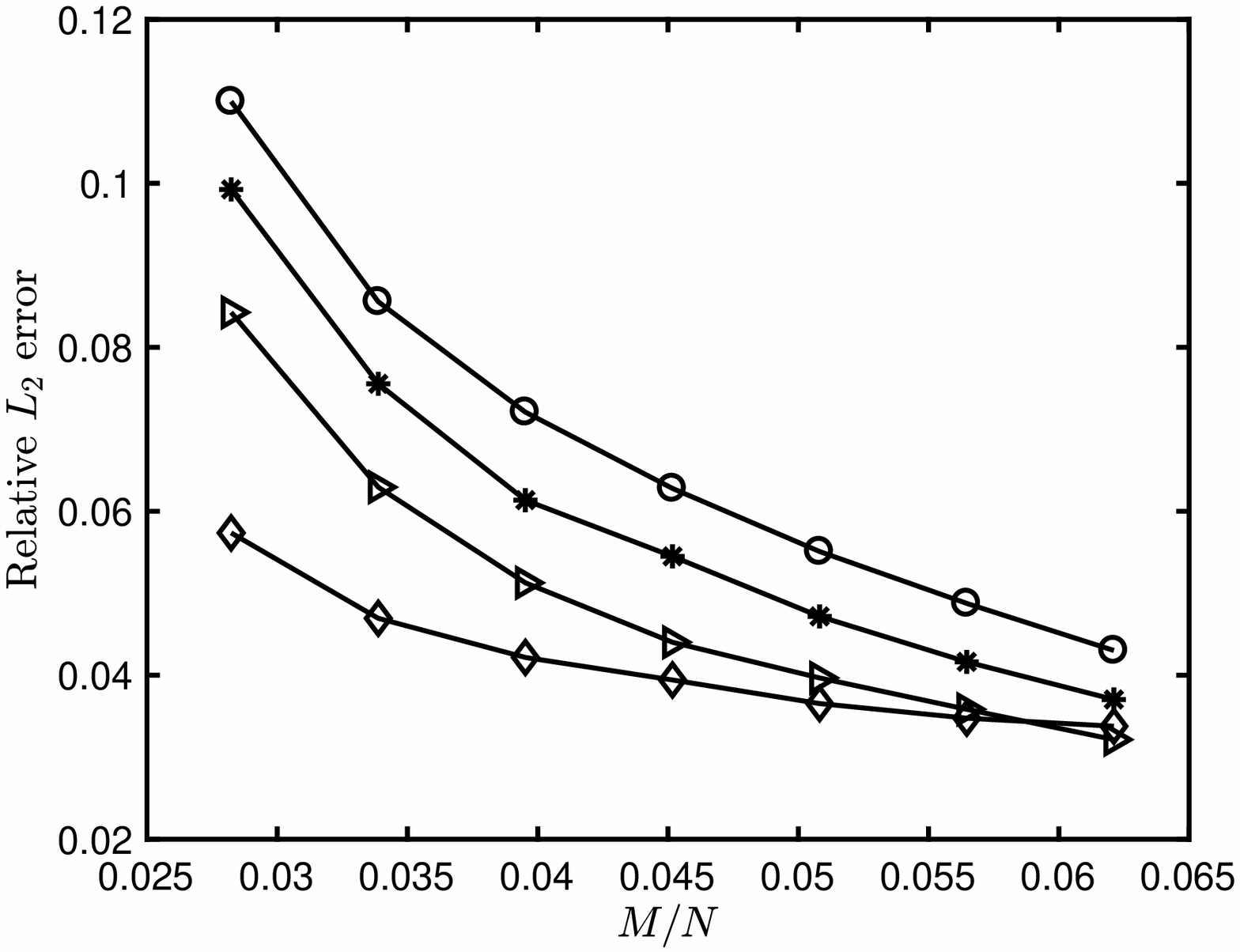}\quad
\includegraphics[width=0.45\textwidth]{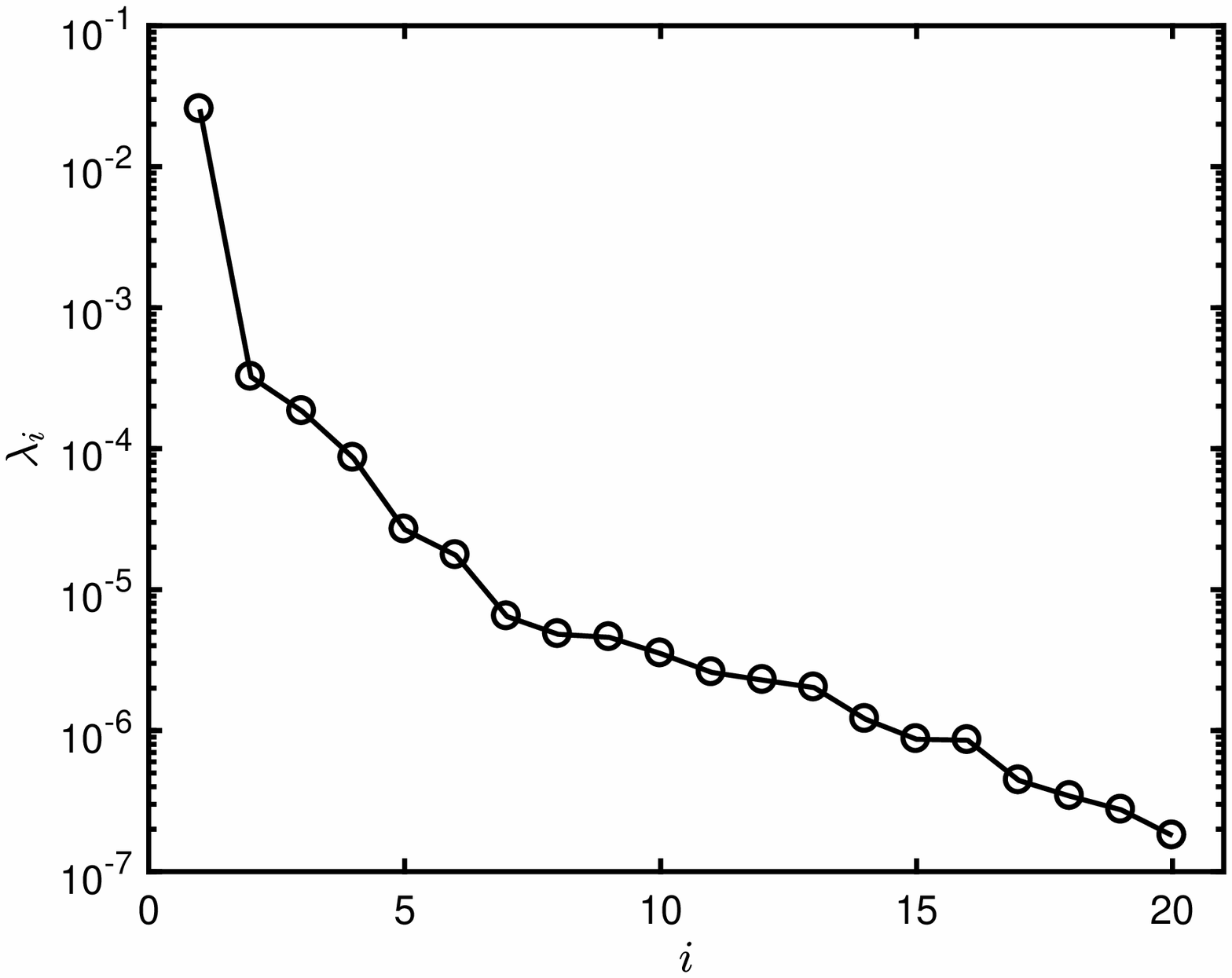}
\caption{Results for the highly compressible function. Left: relative $L_2$
  error. ``$\circ$": standard $\ell_1$, ``$\ast$": re-weighted $\ell_1$, 
  ``$\triangleright$": ADM, ``$\diamond$":
 SADM. Right: eigenvalues of matrix $G$.}
\label{fig:ex2_rmse}
\end{figure}

Consider the following function:
\begin{equation}\label{eq:ex2}
u(\bx) = \sum_{|\ba|=0}^P c_{\ba}\psi_{\ba}(\bx)
       = \sum_{n=1}^Nc_n\psi_n(\bx), \quad
\bx =  (\xi_1,\xi_2,\cdots,\xi_{d}),
\end{equation}
where, $\psi_{\ba}$ are normalized multivariate Hermite polynomials, $d=20, P=3,
N=1771$, and the coefficients $c_n$ are chosen as uniformly distributed random 
numbers,
\begin{equation}
c_n = \zeta/n^{1.5}, \quad \zeta\sim \mathcal{U}[-1,1].
\end{equation}
For this example, we generate $N$ samples of $\zeta$: 
$\zeta^1,\zeta^2,\cdots,\zeta^N$, and then divide them by $n^{1.5}, n=1,2,\cdots,N$
to obtain a random ``compressible signal" $\bm c$. 
This example is also used in \cite{YangLBL16, YangWL17} to demonstrate the
effectiveness of the rotational $\ell_1$ method. The dimension is increased to 
$d=20$ in this test. The function $u$ is not exactly sparse before or after 
rotation. This is reflected in the right plot in Fig.~\ref{fig:ex2_rmse}, which
shows the eigenvalues of $\tensor G$. These eigenvalues indicate that all
subspaces identified by eigen-decomposition of $\tensor G$ make contributions to
$u$, although some of them are quite insignificant. The left plot in 
Fig.~\ref{fig:ex2_rmse} shows results obtained by applying 
Algorithm \ref{algo:cs_rot} and Algorithm \ref{algo:cs_sir1} with re-weighted 
$\ell_1$ minimization and compares them with the standard $\ell_1$ and re-weighted 
$\ell_1$ methods. Apparently, the ADM algorithm improves the accuracy, and SIR 
provides a better initial guess of $\tensor A$, especially when $M/N$ is very 
small, i.e., when the available data are very limited. We also notice that as $M$
increases, the advantage of the SIR-based ADM method decreases.
When $M=110$ ($M/N\approx 0.062$), the SIR-based ADM
(Algorithm~\ref{algo:cs_sir1}) is slightly less accurate than the ADM method 
(Algorithm~\ref{algo:cs_rot}). This is because the size $M$ is sufficient to 
compute a $u_g$ that can yield a slightly better initial guess of $\tensor A$ 
using the gradient information than SIR.

%
\subsection{Korteweg-de Vries equation}
\begin{figure}[h]
\centering
\includegraphics[width=0.45\textwidth]{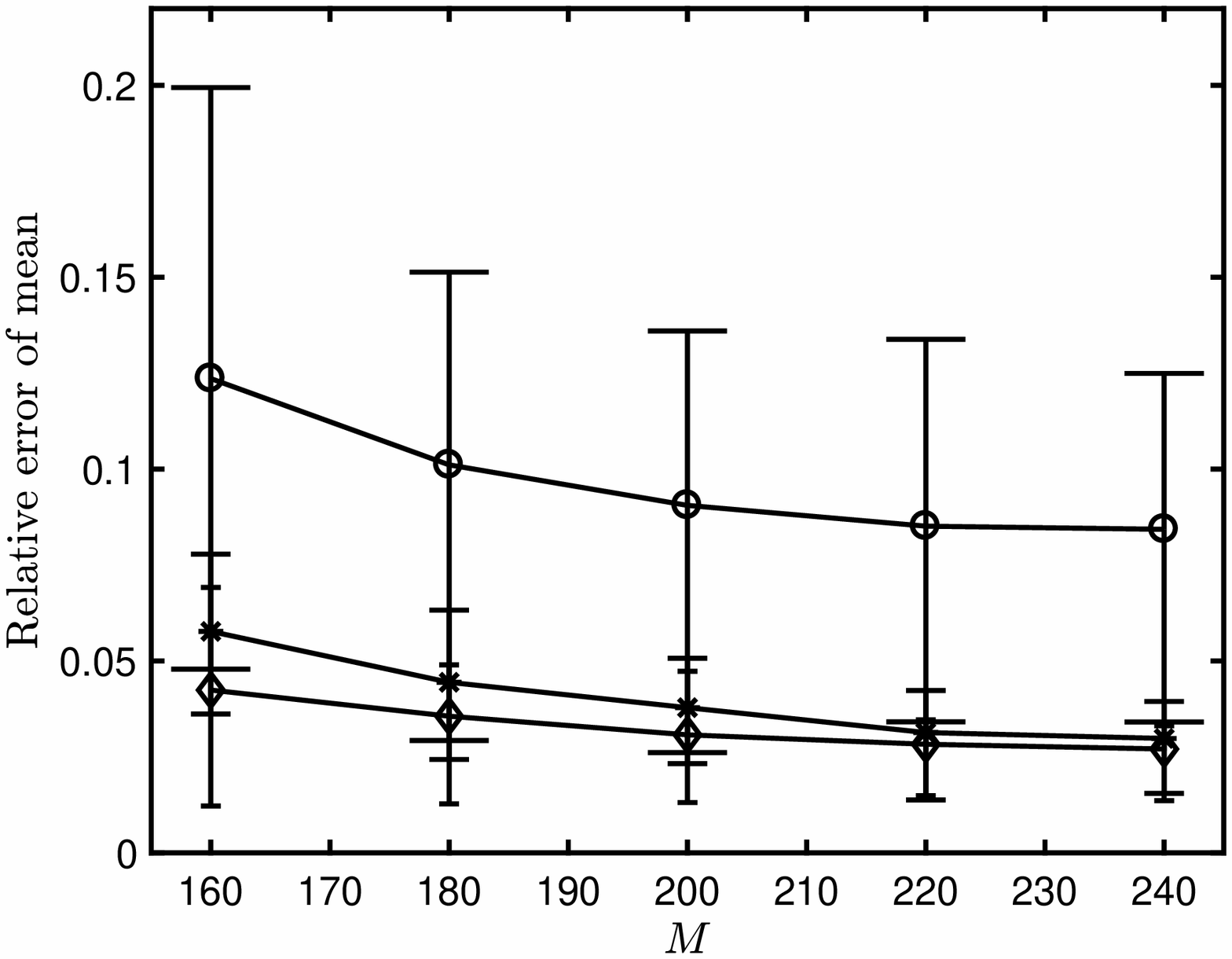}\quad
\includegraphics[width=0.45\textwidth]{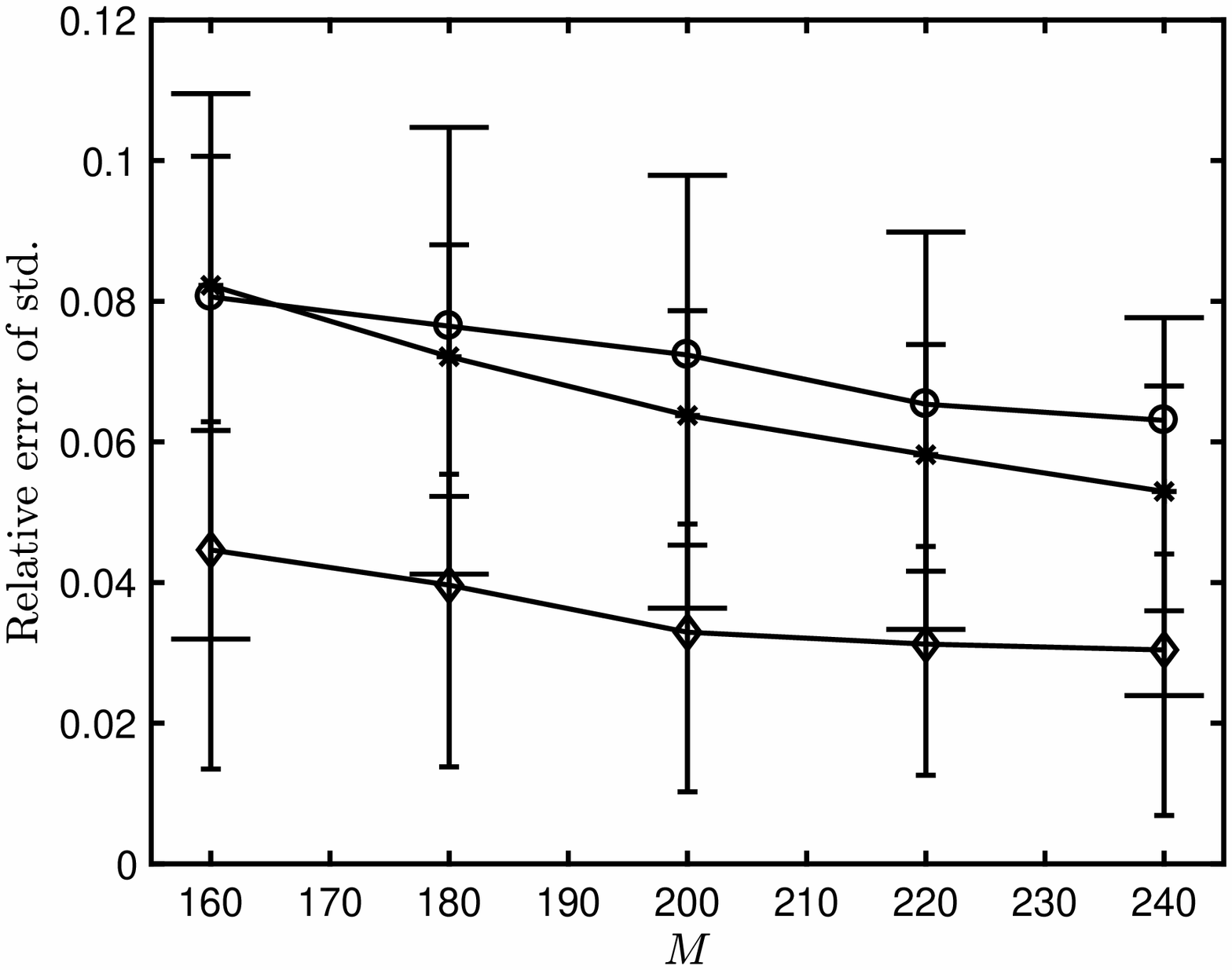}
\caption{Results for the Korteweg-de Vries equation. Left: relative error of
  mean; Right: relative error of standard deviation.
  ``$\circ$": direct estimate from Monte Carlo samples;
``$\ast$": re-weighted $\ell_1$; ``$\diamond$": SADMDR.
``\rule{0.4cm}{0.6pt}": quantiles of direct estimate from Monte Carlo samples;
``\rule{0.2cm}{0.6pt}": quantiles of re-weighted $\ell_1$;
``\rule{0.1cm}{0.6pt}": quantiles of SADMDR.}
\label{fig:ex3_err}
\end{figure}

As an example application of the new method to a nonlinear differential 
equation, we consider the Korteweg-de Vries (KdV) equation with time-dependent 
additive noise \cite{LinGK06}:
\begin{equation}\label{eq:kdv}
\begin{aligned}
& u_t(x,t;\bx)-6u(x,t;\bx)u_x(x,t;\bx)+u_{xxx}(x,t;\bx)=f(t;\bx), 
  \quad x\in (-\infty,\infty), \\
& u(x,0;\bx) = -2 \sech^2(x).
\end{aligned}
\end{equation}
We model $f(t;\bx)$ as a random field represented by the following
Karhuen-Lo\`{e}ve (KL) expansion:
\begin{equation}
f(t;\bx) = \sigma\sum_{i=1}^d\sqrt{\lambda_i}\phi_i(t)\xi_i, 
\end{equation}
where $\sigma$ is a constant and $\{\lambda_i,\phi_i(t)\}_{i=1}^d$ are
eigenpairs of the exponential covariance kernel 
\[C(t, t')=\exp\left(-\dfrac{|t-t'|}{l_c}\right).\]
The explicit form of $\lambda_i$ and $\phi_i$ can be found in \cite{GhanemS91}. 
In this problem, we set $l_c=0.1$ and $d=100$ 
($\sum_{i=1}^d\lambda_i > 0.978\sum_{i=1}^{\infty}\lambda_i$). In this 
case, the exact one-soliton solution is
\begin{align}\label{eq:kdv_sol}
  \begin{split}
    u(x,t;\bx)=&\sigma\sum_{i=1}^d\sqrt{\lambda_i}\xi_i\int_0^t\phi_i(y)\dif y\\
    & - 2\sech^2\left(x-4t+6\sigma\sum_{i=1}^d\sqrt{\lambda_i}\xi_i
\int_0^t\int_0^{z}\phi_i(y)\dif y\dif z\right).
\end{split}
\end{align}
The QoI is chosen to be $u(x,t;\bx)$ at $x=6,t=1$ with 
$\sigma=0.4$. Because an analytical expression for $\phi_i$ is available, 
we can compute the integrals in Eq.~\eqref{eq:kdv_sol} with high accuracy. Denoting
\begin{equation}
A_i = \sqrt{\lambda_i}\int_0^1\phi_i(y)\dif y,\quad 
B_i = \sqrt{\lambda_i}\int_0^1\int_0^{z}\phi_i(y)\dif y\dif z, \quad
i=1,2,\cdots,d,
\end{equation}
the analytical solution is
\begin{equation}\label{eq:kdv_sol2}
u(x,t;\bx)\big |_{x=6,t=1}=\sigma\sum_{i=1}^dA_i\xi_i
-2\sech^2\left(2+6\sigma\sum_{i=1}^d B_i\xi_i\right).
\end{equation}
We set $P=2$ ($N=5151$) to construct $u_g$ using re-weighted $\ell_1$ 
minimization and set $\tilde d=12$ and $P=3$ ($N=455$) to construct $\tilde u_g$
using SADMDR (Algorithm~\ref{algo:cs_sir2}).
The relative error of the mean and standard deviation obtained from ADM, SADMDR,
and the MC method with $M$ realizations, compared with the reference solution, 
are presented in Fig.~\ref{fig:ex3_err}. For the estimate of mean, both 
re-weighted $\ell_1$ and SADMDR are more accurate than MC, and the 
SADMDR is up to $30\%$ more accurate than re-weighted $\ell_1$ for small $M$. 
The accuracy of ADM and SADMDR becomes
similar as $M$ increases. For the estimate of the standard deviation, the advantage
of SADMDR over ADM is much more distinct: for all considered $M$, SADMDR has a
$50\%$ smaller error than MC, while the re-weighted $\ell_1$ method has a similar 
error as MC for small $M$ (error of re-weighted $\ell_1$ is slightly larger than 
in MC for $M=160$) and $20\%$ smaller error than MC for larger $M$. Again, the
observed difference between the re-weighted $\ell_1$ and SADMDR results 
become smaller as $M$ increases.
Moreover, the results obtained from Algorithm~\ref{algo:cs_rot} and
Algorithm~\ref{algo:cs_sir1} are almost the same as the re-weighted $\ell_1$
results, 
so the results are not plotted. This implies that
the rotations identified in the iteration by the gradient of $u$ approximated 
from $u_g$ with up to second-order Hermite polynomials are not sufficiently
informative to provide good guidance for a sparser representation--no matter what
initial guess is used. In \cite{YangLBL16}, the KdV equation with 
$12$ random parameters is solved using a $u_g$ approximation with $P=4$. In
that example, Algorithm~\ref{algo:cs_rot} reduced the relative $L_2$ error of 
re-weighted $\ell_1$ by up to $75\%$. This implies that higher-order terms in 
the Hermite polynomial expansion are important for determining the rotation 
matrix. In the KdV equation used in this work, $d=100$, and $M$ is small. Thus,
there are not enough samples to compute terms even in $P=3$ gPC expansion 
without first reducing the dimensionality $d$. We use SIR to set reduced 
dimension $\tilde{d}=12$. Then, we choose $P=3$ to include as many terms as 
possible in the Hermite polynomial expansion but still keep the number of 
unknowns in a reasonable range ($N=455$).

%
\subsection{Groundwater flow}
\label{subsec:ex4}
\begin{figure}[h]
\centering
\includegraphics[width=0.46\textwidth]{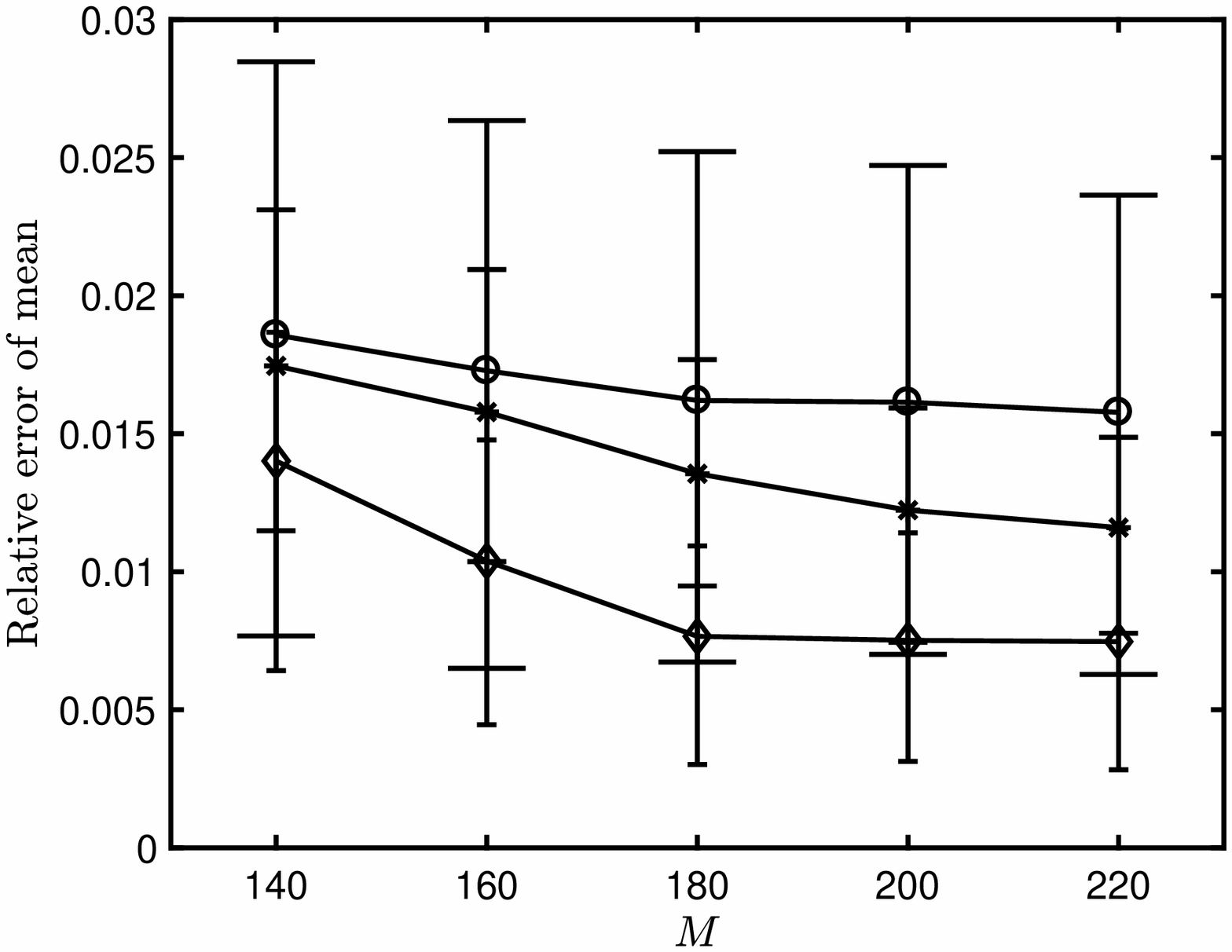}\quad
\includegraphics[width=0.45\textwidth]{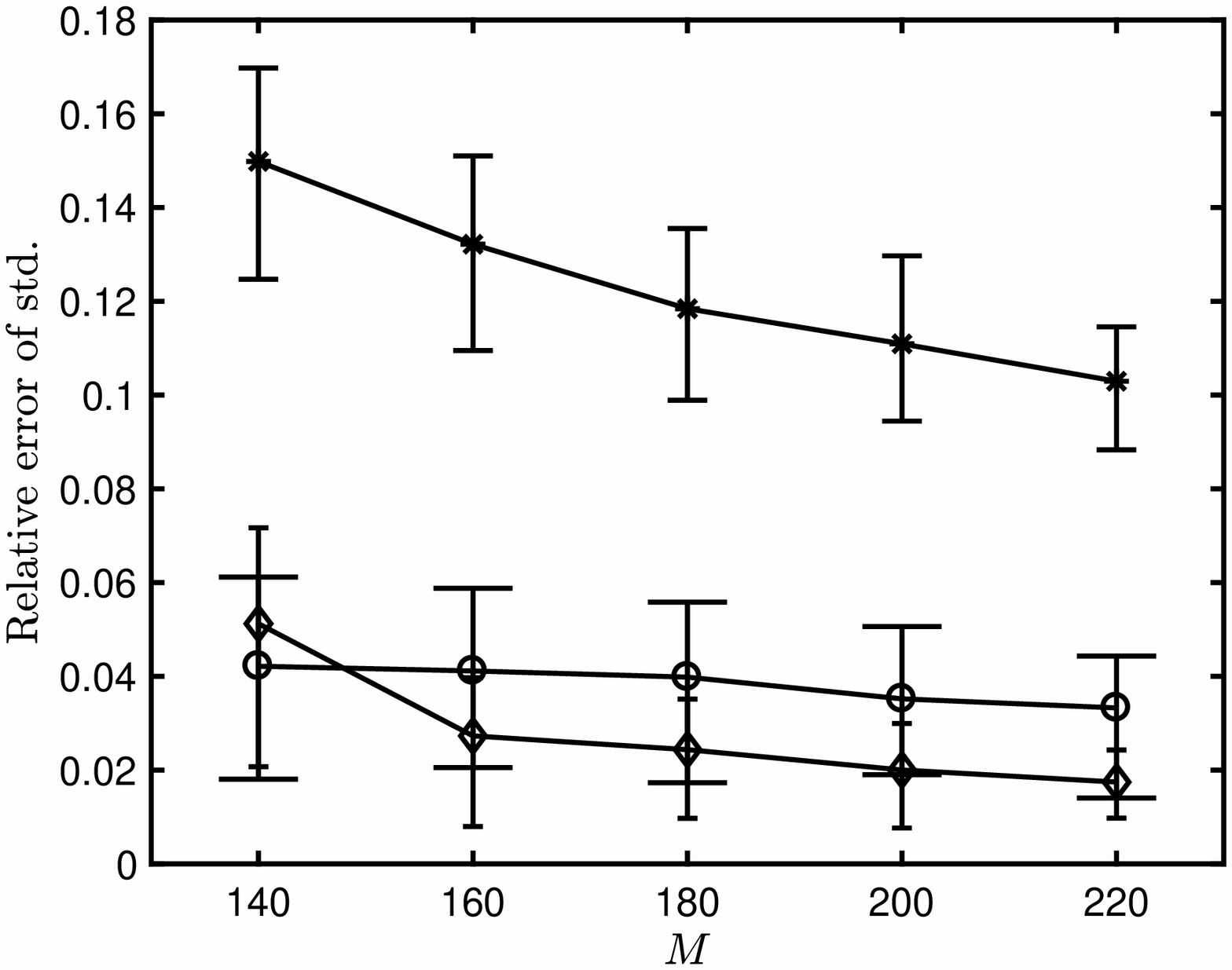}
\caption{Results for the groundwater flow. Left: relative error of
  mean; Right: relative error of standard deviation. ``$\circ$": direct estimate 
from Monte Carlo samples; ``$\ast$": re-weighted $\ell_1$; ``$\diamond$": SADMDR.
``\rule{0.4cm}{0.6pt}": quantiles of direct estimate from Monte Carlo samples;
``\rule{0.2cm}{0.6pt}": quantiles of re-weighted $\ell_1$;
``\rule{0.1cm}{0.6pt}": quantiles of SADMDR.}
\label{fig:ex4_error}
\end{figure}

Next, we consider a model that simulates the groundwater flow in a confined
aquifer, which spans a $2000$ m $\times 1000$ m area \cite{LiLinLi16}. The north 
and south boundaries of the aquifer are two rivers with constant but different 
hydraulic heads, whereas the east and west boundaries are bounded by no-flow 
conditions. This model can be described by the following equations:
\begin{equation}
  \begin{cases}
    q(x,y) = -T(x,y)\nabla u(x,y), \quad (x,y)\in D=[0,2000]\times [0, 1000]\\
    \nabla\cdot q(x,y) = 0, 
  \end{cases}
\end{equation}
and boundary conditions
\begin{equation}
  \begin{cases}
    u(x,0) = 0, ~ u(x, 1000) = 10, \\
    q_x(0,y) = q_x(2000, y) = 0,
  \end{cases}
\end{equation}
where $u$ is hydraulic head $[\text{m}]$, $T$ is transmissivity $[\text{m}^2/\text{day}]$,
and $\bm q=(q_x, q_y)^\trans$ is flux vector $[m^2/\text{day}]$. The QoI 
is the hydraulic head at a specific location: $u^*=u(200, 500)$. The equations
are solved by the finite difference method on a $61\times 31$ computational
grid. The transmissivity $T(x,y)$ is described by log-normal random field 
$T(x,y)=\ln S(x,y)$, where $S$ satisfies: 1) for $(x,y)\in D$, $S(x,y)\sim
\mathcal{N}(2, 1)$; 2) for $(x,y), (x',y')\in D$, the covariance kernel is
\[C(x,y;x',y')=\exp\left(-\dfrac{|x-x'|}{l_x}-\dfrac{|y-y'|}{l_y}\right).\]
In this problem, we set $l_x=l_y=300$, and $S(x,y)$ is represented by a KL-expansion 
with $100$ terms ($\sum_{i=1}^{100}\lambda_i > 0.85\sum_{i=1}^{\infty}\lambda_i$).
We set $P=2$ ($N=5151$) to construct $u_g$ using the re-weighted $\ell_1$
method and set $\tilde d=20$ and $P=3$ ($N=1771$) to construct $\tilde u_g$ 
using SADMDR (Algorithm~\ref{algo:cs_sir2}). Figure~\ref{fig:ex4_error}
represents the relative error of mean and standard deviation estimates. For the 
mean, both methods exhibit better accuracy than direct estimation from MC
with the same number of sampling points. The re-weighted $\ell_1$ method reduces 
the relative error by up to $30\%$ compared with MC, and the SADMDR reduces the
relative error by up to $50\%$. For the standard deviation, re-weighted $\ell_1$
has relative error that is several times larger than the error of MC for all 
considered $M$. On the other hand, SADMDR reduces the error of MC by 
approximately $45\%$ for $M>160$. For $M=140$, SADMDR performs worse than MC 
because the sample size is too small for this high-dimensional ($d=100$) problem.
Again, the Algorithms~\ref{algo:cs_rot} and \ref{algo:cs_sir1} results are 
not presented as they are almost the same as those by re-weighted $\ell_1$. This
example also demonstrates that reducing dimension and increasing $P$ while 
keeping $N$ in a reasonable range for compressive sensing is important for an 
accurate approximation of $u$.


\subsection{High-dimensional function}
\begin{figure}[h]
\centering
\includegraphics[width=0.45\textwidth]{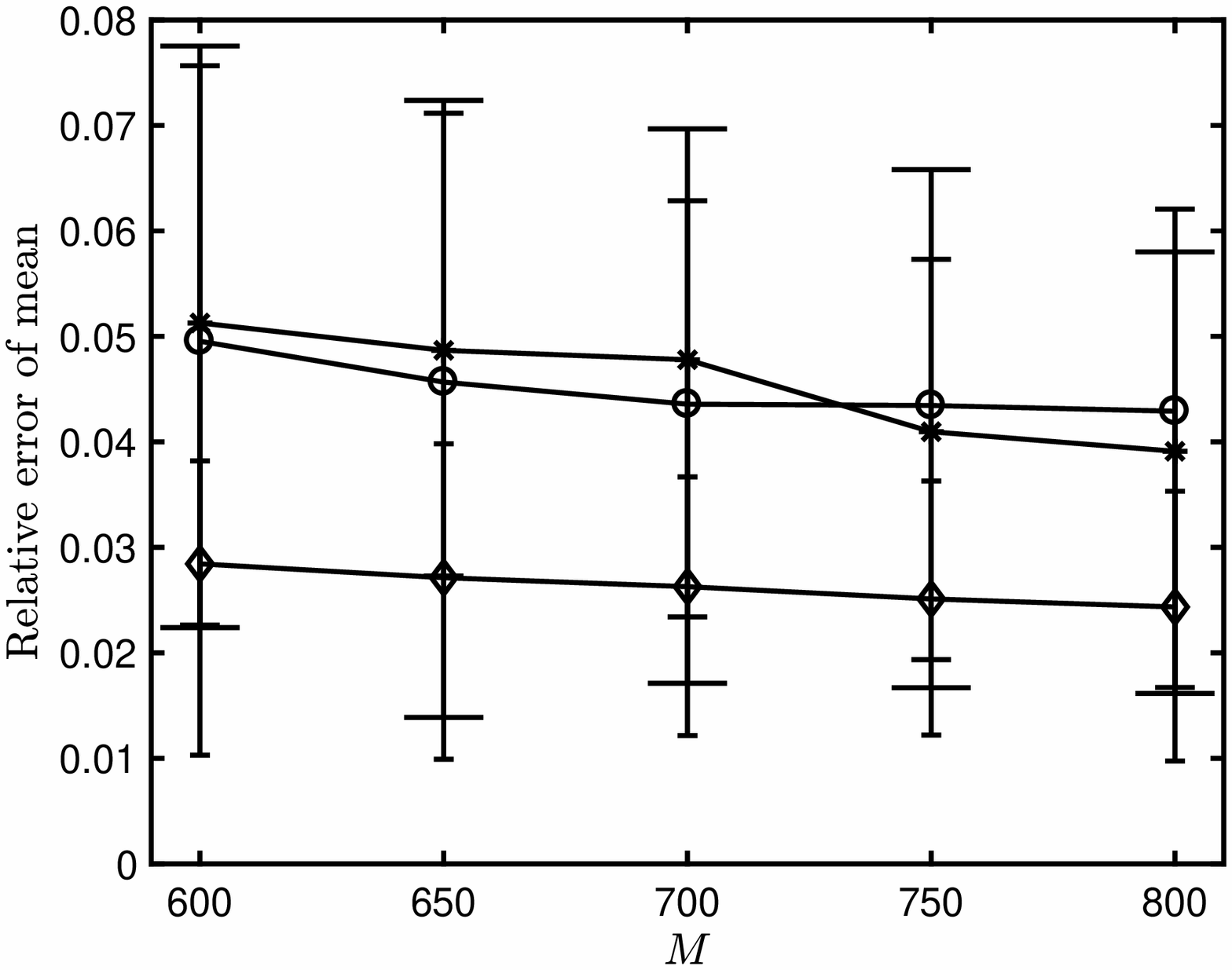}\quad
\includegraphics[width=0.45\textwidth]{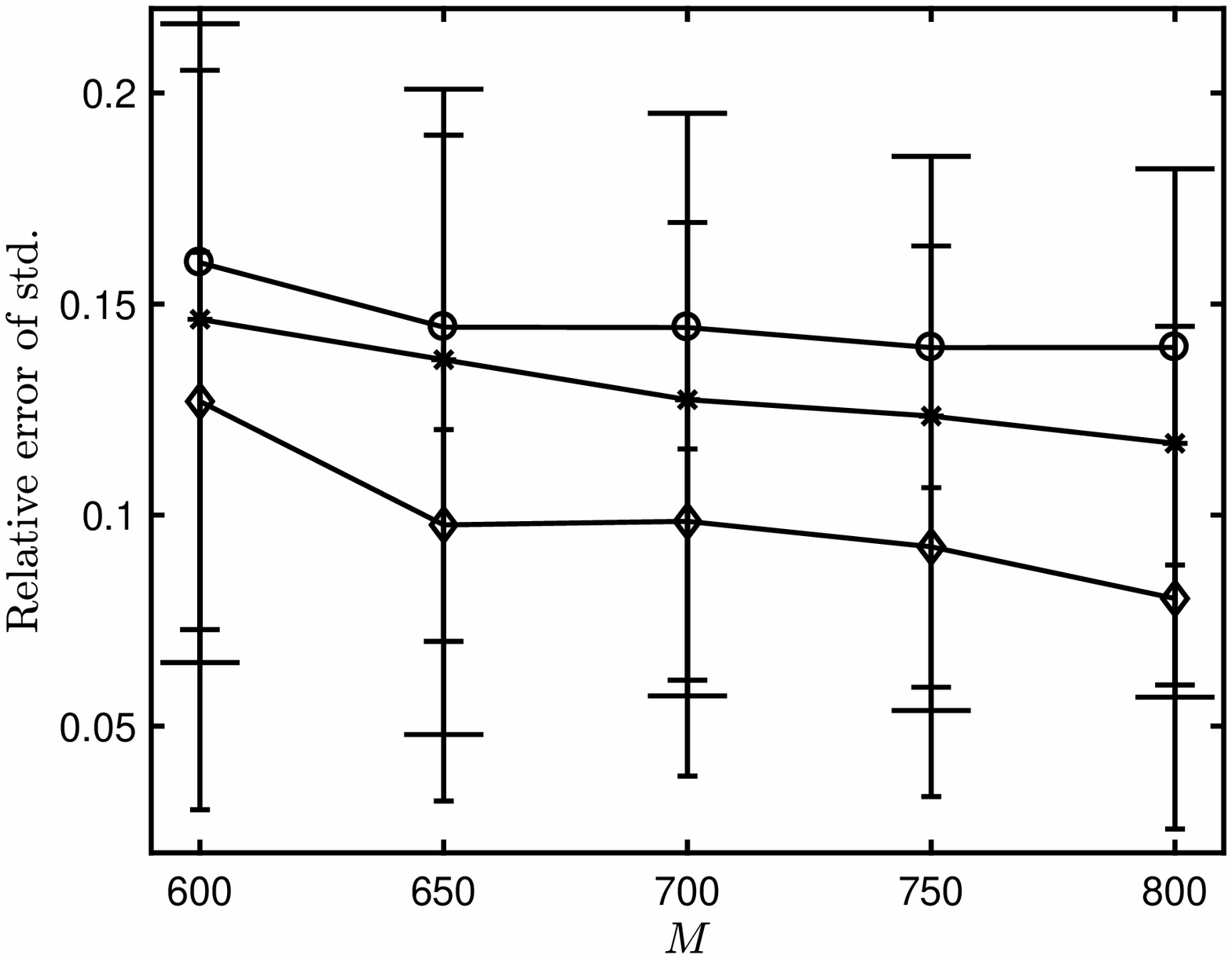}
\caption{Results for the high-dimensional function. ``$\circ$": direct estimate 
from Monte Carlo samples; ``$\ast$": re-weighted $\ell_1$; ``$\diamond$": SADMDR.
``\rule{0.4cm}{0.6pt}": quantiles of direct estimate from Monte Carlo samples;
``\rule{0.2cm}{0.6pt}": quantiles of re-weighted $\ell_1$;
``\rule{0.1cm}{0.6pt}": quantiles of SADMDR.}
\label{fig:ex5_error}
\end{figure}

In the final example, we demonstrate the ability of SADMDR to deal with very 
high-dimensional problems. Specifically, we consider the
following function \cite{hampton2018basis}: 
\begin{equation}
  u(\bx) = \exp\left(2-\sum_{i=1}^d\dfrac{\sin(i)\xi_i}{i}\right), \quad d=500.
\end{equation}
We use a third-order gPC expansion without interaction terms, i.e.,  
we only use constant and $\{\xi_i, (\xi_i^2-1)/\sqrt{2}, (\xi_i^3-3\xi_i)/\sqrt{6}\}_{i=1}^d$
as basis functions ($N=1+500+500=1001$), to construct $u_g$ using the re-weighted $\ell_1$ 
method, then we set $\tilde d=20$ and $P=3$ ($N=1771$) to construct $\tilde u_g$ using 
SADMDR (Algorithm~\ref{algo:cs_sir2}). The relative errors of the mean and 
standard deviation are presented in Fig.~\ref{fig:ex5_error}. As before, SADMDR 
reduces the error in the mean prediction by $40\%$ compared with MC.
The $u_g$ by re-weighted $\ell_1$ provide a similar error in the estimate of
mean as MC. For the standard deviation estimate, the error in re-weighted $\ell_1$ 
is approximately $10\%$ smaller than $MC$, while SADMDR has approximately $30\%$
smaller error than in MC. More importantly, unlike previous examples where the 
differences between re-weighted $\ell_1$ and SADMDR became smaller quickly with 
increasing $M$, in this case, the accuracy of SADMDR relative to re-weighted $\ell_1$ 
changes slowly in the range of studied $M$. This is because the dimension of the
problem is very high, and dimension reduction is more critical than in the 
previous examples.
The ADM and SADM algorithms fail to improve the accuracy in this example as in
Examples 3 and 4, because the dimension is very high and the available data are
too limited for these two approaches.

\subsection{Discussion}
\begin{figure}[h]
\centering
\includegraphics[width=0.45\textwidth]{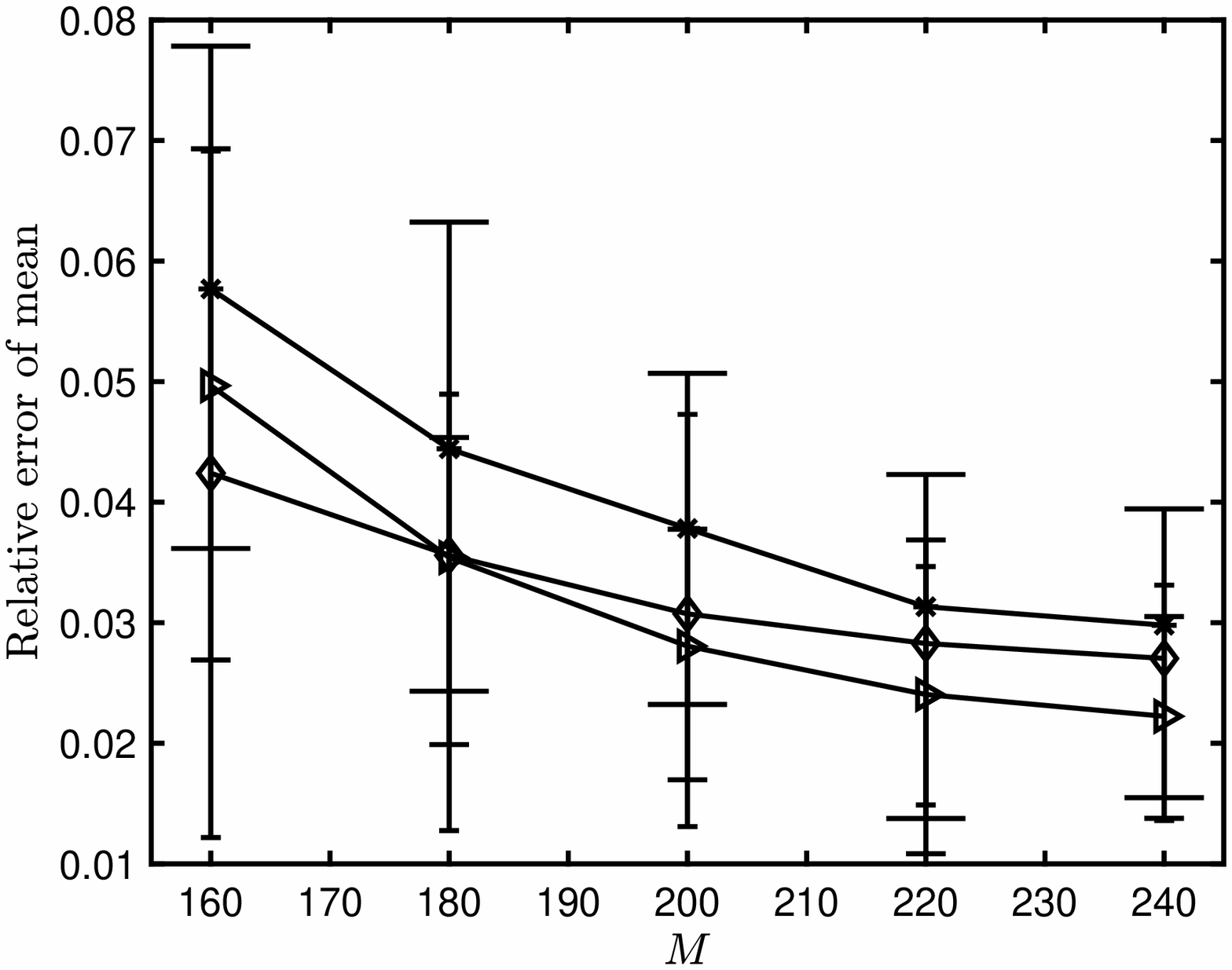}\quad
\includegraphics[width=0.45\textwidth]{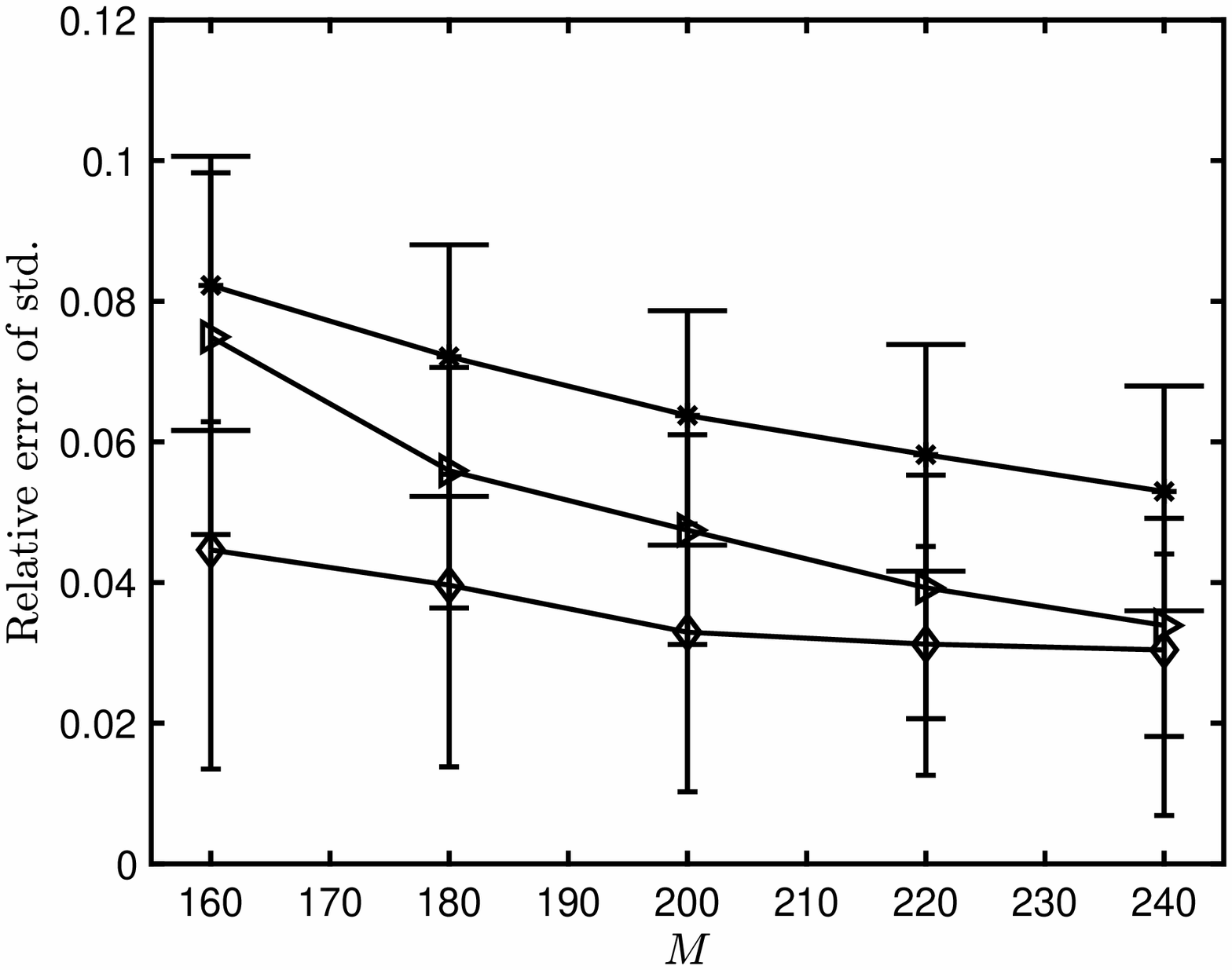}
\caption{Results for the high-dimensional function.  
``$\ast$": re-weighted $\ell_1$; ``$\rhd$" dimension reduction based on re-weight
$\ell_1$;``$\diamond$": SADMDR.
``\rule{0.4cm}{0.6pt}": quantiles of re-weighted $\ell_1$ ($u_g$);
``\rule{0.2cm}{0.6pt}": quantiles of gradient-based dimension reduction ($\tilde u_g^W$);
``\rule{0.1cm}{0.6pt}": quantiles of SADMDR ($\tilde u_g$).}
\label{fig:ex6}
\end{figure}

The rotation matrix computed from $\tensor G$ in Eq.~\eqref{eq:grad_mat} was
used for dimension reduction in active subspace method \cite{Russi10, ConstantineDW14}. 
We can also use it to replace SIR in Algorithm \ref{algo:cs_sir2}. 
Specifically, on computing the eigen-decomposition of $\tensor G$ based on the $u_g$ from
Algorithms~\ref{algo:cs1} or \ref{algo:cs_rot}, we can construct $\hat{\tensor A}$ that consists
of eigenvectors corresponding to the $\tilde d$ largest eigenvalues, 
where $\tilde d < d$ is the reduced dimension. Then we project $\bm\xi$ to 
$\tilde{\bm\xi} = \hat{\tensor A}\bm\xi$, and run step 4\ in 
Algorithm~\ref{algo:cs_sir2} to construct a surrogate model. Apparently, the
accuracy of estimating $\nabla u$ is critical to the dimension reduction. Since
we do not have samples of $\nabla u$, we have
to approximate $\nabla u$ from $\nabla u_g$. Therefore one can roughly expect 
that an accurate $u_g$ yields better performance of dimension reduction, while a
less accurate $u_g$ results in worse performance. We use the KdV equation in
Example 3 as a demonstration. As we presented in Example 3, we set $P=2$ 
($N=5151$) to construct $u_g$ using re-weighted $\ell_1$ minimization. Then
we approximate $\tensor G$ based on $u_g$, and truncate the dimension to 
$\tilde d=12$. Next, we set $P=3$
($N=455$), and construct $\tilde u_g^W$ using
ADM method. We compare the accuracy of this $\tilde u_g^W$ with $u_g$ and
$\tilde u_g$ (from SADMDR by setting $\tilde d=12$) in Example 3 and present the
results in Fig.~\ref{fig:ex6}. When $M<200$, $\tilde u_g$ is more accurate in 
estimating mean than $u_g$ and $u_g^W$. But when $M\geq 200$, $\tilde u_g^W$ is
the best for estimating the mean. For the estimate of the standard deviation, 
$\tilde u_g$ is always the best in the range of $M$ we chose. But the difference
between $\tilde u_g$ and $\tilde u_g^W$ decays as $M$ increases. These phenomena 
are similar to those in Example 2. Again, this is because as $M$ increases, 
$\ell_1$ or re-weighted $\ell_1$ minimization is able to provide more accurate 
estimate of $u_g$, and consequently more accurate estimates of $\nabla u$ and 
$\tensor G$. Therefore, the initial guess of $\tensor A$ (Example 2) or 
$\hat{\tensor A}$ (this example) becomes better, and finally exceeds the one 
estimated from SIR when $M$ is sufficiently large. 
An approach using Algorithm~\ref{algo:cs_rot} to estimate $\tensor G$, then 
perform dimension reduction was proposed in \cite{YangBRT17}. It worked well
for specific problems in that study. As we show in Example 2 and this
comparison, there is no guarantee that SIR works better than the 
gradient-based method for dimension reduction for any $M$. If the gradient 
information is not available, or it is difficult to approximate the gradient 
accurately, SIR can be a good choice. In practice, if no prior knowledge of the 
QoI is available, users may construct different surrogate models, then use
model selection tools (e.g., AIC \cite{Akaike74}, BIC \cite{Schwarz1978} or cross
validation \cite{Kohavi1995}) to decide which method to employ.


\section{Conclusions}
\label{sec:conclusion}

We use the sliced inverse regression method to provide a better initial
guess for the alternating direction method proposed in \cite{YangLBL16}, which
enhances the sparsity of the Hermite polynomial expansion of the QoI relying on 
i.i.d. Gaussian random variables. The enhancement of sparsity helps the 
compressive sensing method to obtain a more accurate Hermite polynomial 
expansion. Examples 1 and 2 show that when the available data (e.g., the number
of model realizations, $M$) are limited (compared with the number of unknown 
terms in the QoI expansion), the compressive sensing method with the initial 
guess provided by SIR can yield more accurate results than the standard $\ell_1$ 
minimization (or re-weighted $\ell_1$) used in \cite{YangLBL16}. We also 
demonstrate that when the problem dimensionality is very high and the available
data size is small, we can first use SIR to perform dimension reduction then use
the ADM method based on the reduced system to approximate the mean and standard 
deviation of the QoI more accurately. The dimension reduction allows for
inclusion of more higher-order terms in the Hermite polynomial expansion of the
QoI. Consequently, more information related to variance and other higher-order 
terms are included in the expansion. This is illustrated in Examples 3-5 where 
the improvement in the estimate of standard deviation is much more significant 
than in the estimate of the mean. We demonstrate the advantage of the new 
algorithms (i.e., Algorithms~\ref{algo:cs_sir1} and \ref{algo:cs_sir2})  
when the number of samples is much smaller than the number of terms in the QoI
expansion and is far below the requirement of the sample size for the $\ell_1$
minimization. Examples 1 and 2 (relatively low-dimensional problems) demonstrate
that the accuracy of the ADM method (i.e., Algorithm~\ref{algo:cs_rot}) improves
with increasing $M$ and it can be as good as or even better than
the SIR-aided method for larger $M$. In Examples 3-5 (higher-dimensional problems), $\ell_1$ 
minimization without dimension reduction does not work well for relatively small $M$.
However, if we dramatically increase $M$, it is expected that the accuracy of 
ADM will be comparable to that of SIR-aided ADM methods.

In this work, the rotation matrix is obtained in two different ways. SIR uses
conditional mean to identify the matrix $\hat A$, while in ADM, the gradient of
$u$ is used to iteratively identify the rotation. Both approaches have been 
widely used in statistics algorithms for dimension reduction. Notably, other 
approaches also can be incorporated in our algorithm. For example, 
different methods for sufficient dimension reduction 
(e.g., \cite{Li2007sparse, JiangL14}) may provide a better initial guess of the
rotation matrix or better dimension reduction strategy for specific problems.
In each iteration, the rotation can be obtained from these SIR-type methods 
instead of using the gradient information. Also, when $d$ is not very large 
(typically $d<100$), Algorithm \ref{algo:cs_rot} or Algorithm \ref{algo:cs_sir1}
can be used to obtain $u_g$ then construct the gradient matrix $\tensor G$ of
$u$ and reduce the dimension according to the magnitude of eigenvalues of 
$\tensor G$ \cite{YangBRT17}. 

Moreover, we demonstrate the effectiveness of ADM for $\ell_1$ minimization.
ADM can also be integrated with other optimization methods to solve
the compressive sensing problem, e.g., OMP \cite{BrucksteinDE09}, $\ell_{1-2}$
minimization \cite{YinLHX15}, etc. Further, it could be advantageous to integrate
our method with sampling strategies (e.g., \cite{AlemazkoorM17, AlemazkoorM17-2}),
basis selection method (e.g., \cite{JakemanES14}), or Bayesian approach 
(e.g., \cite{KaragiannisBL15}). The combination of these methods can be
especially useful for problems where experiments or simulations are costly and
where a good surrogate model of the QoI is needed, e.g., in inverse problems based
on a Bayesian framework (\cite{SargsyanSNDRT14,YangLGTMB17}). 

Finally, as discussed in Section 3, a correct balance between the reduced
dimension $\tilde d$ and the selection of high-order terms in the expansion can
yield a more accurate approximation of the QoI. The theoretical analysis on 
SIR-type dimension reduction methods can be found 
in \cite{Li2007sparse, JiangL14}, which help to identify $\tilde d$. 
After $\tilde d$ is set, model selection techniques 
can be used to select the polynomial order $P$. 
In practice, whether to use dimension reduction depends on the
  available data size, model complexity and property of QoI. Again, model
  selection tools can be used to identify a suitable surrogate model when no
  prior knowledge is available.

\section*{Appendix}
\renewcommand{\theequation}{A-\arabic{equation}}

\appendix

\subsection*{A. Cross-validation method}

The algorithm in \cite{DoostanO11} is used to estimate the error term $\epsilon$ 
in $(P_{1,\epsilon})$. This algorithm is summarized in Algorithm~\ref{algo:cross}.
\begin{algorithm}
\caption{Cross-validation to estimate the error $\epsilon$.}
\label{algo:cross}
\begin{algorithmic}[1]
\STATE Divide the $M$ output samples to $M_r$ reconstruction ($\bm u_r$) and
$M_v$ validation ($\bm u_v$) samples and divide the measurement matrix 
$\tensor\Psi$ correspondingly into $\tensor\Psi_r$ and $\tensor\Psi_v$.
\STATE Choose multiple values for $\epsilon_r$ such that the exact error
$\Vert\tensor\Psi_r\bm c-\bm u_r\Vert_2$ of the reconstruction samples is
within the range of $\epsilon_r$ values.
\STATE For each $\epsilon_r$, solve $(P_{1,\epsilon})$ with $\bm u_r$ and
$\tensor\Psi_r$ to obtain $\hat{\bm c}$. Then compute 
$\epsilon_v=\Vert\tensor\Psi_v\hat{\bm c}-\bm u_v\Vert_2$.
\STATE Find the minimum value of $\epsilon_v$ and its corresponding 
$\epsilon_r$. Set $\epsilon=\sqrt{M/M_r}\epsilon_r$.
\end{algorithmic}
\end{algorithm}
Of note, a technique to avoid the cross-validation step in some cases is
proposed in \cite{Adcock17},. 


\section*{Ackowledgement}
This work was supported by the U.S. Department of Energy (DOE), Office of 
Science, Office of Advanced Scientific Computing Research (ASCR) as part of the
Multifaceted Mathematics for Complex Systems project and the Uncertainty 
Quantification in Advection-Diffusion-Reaction Systems projects. A portion of 
the research described in this paper was conducted under the Laboratory Directed
Research and Development Program at Pacific Northwest National Laboratory 
(PNNL). PNNL is operated by Battelle for the DOE under Contract DE-AC05-76RL01830.

\bibliographystyle{plain}
\bibliography{uq}

\end{document}